\documentclass[oneside,11pt,dvipsnames]{article}

\usepackage{arxiv}

\usepackage[T1]{fontenc}      

\usepackage{lmodern}
\usepackage{tikz}           
\usepackage[]{xcolor}       
\usepackage{amsmath}        
\usepackage{amsfonts}       
\usepackage{dsfont}         
\usepackage{hyperref}       
\usepackage{url}            
\usepackage{booktabs}       
\usepackage{nicefrac}       
\usepackage{microtype}      
\usepackage[frozencache]{minted}         
\usepackage{graphicx}       
\usepackage{doi}            
\usepackage{bm}             
\usepackage{printlen}
\usepackage{authblk}
\usepackage[symbol]{footmisc}
\usepackage[toc,page,titletoc]{appendix}

\usepackage[style=numeric-comp,sorting=none,giveninits=true,maxbibnames=10,backend=biber]{biblatex}
\AtEveryBibitem{
    \iffieldundef{doi}{}{\clearfield{url}} 
}
\DeclareCiteCommand{\citenum}
  {\usebibmacro{cite:init}%
   \iffieldundef{prenote}
     {}
     {\BibliographyWarning{Ignoring prenote argument}}%
   \iffieldundef{postnote}
     {}
     {\BibliographyWarning{Ignoring postnote argument}}}%
  {\usebibmacro{citeindex}%
   \usebibmacro{cite:comp}}
  {}
  {\usebibmacro{cite:dump}}

\renewcommand{\thefootnote}{\fnsymbol{footnote}}
\renewcommand{\headeright}{}  
\renewcommand{\undertitle}{}  

\usetikzlibrary{calc,decorations.pathmorphing,decorations.pathreplacing,fadings,shadings,shapes,arrows,positioning,patterns,automata,shadows.blur,fit,arrows.meta} 

\tikzset{
	>=stealth',
	punkt/.style={
		rectangle,
		rounded corners,
		draw=black, very thick,
		text width=12em,
		minimum height=2em,
		text centered,
		fill=black!5},
	pil/.style={
		->,
		thick,
		shorten <=2pt,
		shorten >=2pt,}
}
\tikzstyle{class}=[
rectangle,
draw=black,
text centered,
anchor=north,
text=black,
text width=7cm,
shading=axis,
blur shadow={shadow blur steps=5},
bottom color=black!10,top color=white,shading angle=45]
\tikzstyle{whitebox}=[rectangle, draw=black, text centered, anchor=north, text=black, text width=7cm, bottom color=white,top color=white]
\tikzstyle{estimator} =[class, bottom color=blue!10, top color=blue!10, shading angle=45]
\tikzstyle{class2} = [class, text width=4.7cm]
\tikzstyle{class3} = [class, text width=3.2cm]
\tikzstyle{inharrow}=[->, >=open triangle 90, thick]
\tikzstyle{comp}=[->, >=diamond, thick]

\DeclareMathOperator{\sech}{sech}
\DeclareMathOperator{\law}{law}
\DeclareMathOperator{\diag}{diag}
\DeclareMathOperator{\argmin}{argmin}

\newcommand\extrafootertext[1]{%
	\bgroup
	\renewcommand\thefootnote{\fnsymbol{footnote}}%
	\renewcommand\thempfootnote{\fnsymbol{mpfootnote}}%
	\footnotetext[\value{footnote}]{#1}%
	\egroup
}

\newminted{python}{frame=lines,fontsize=\footnotesize,framerule=.5pt,mathescape,autogobble,tabsize=4,obeytabs}

\author[1]{Moritz Hoffmann}
\author[1]{Martin Scherer}
\author[1,2]{Tim Hempel}
\author[1]{Andreas Mardt}
\author[3,\footnote{}]{Brian de Silva}
\author[1,4,5,6]{Brooke E. Husic}
\author[7]{Stefan Klus}
\author[8]{Hao Wu}
\author[3]{Nathan Kutz}
\author[9]{Steven L. Brunton}
\author[1,2,10]{Frank Noé}

\affil[1]{Fachbereich Mathematik und Informatik, Freie Universität Berlin, 14195 Berlin, Germany}
\affil[2]{Fachbereich Physik, Freie Universität Berlin, 14195 Berlin, Germany}
\affil[3]{Department of Applied Mathematics, University of Washington, Seattle, WA 98105, United States}
\affil[4]{Lewis Sigler Institute for Integrative Genomics, Princeton University, Princeton, NJ 08540, United States}
\affil[5]{Princeton Center for Theoretical Science, Princeton University, Princeton, NJ 08540, United States}
\affil[6]{Center for the Physics of Biological Function, Princeton University, Princeton, NJ 08540, United States}
\affil[7]{Department of Mathematics, University of Surrey, Guildford, GU2 7XH, United Kingdom}
\affil[8]{School of Mathematical Sciences, Tongji University, Shanghai 200092, P. R. China}
\affil[9]{Department of Mechanical Engineering, University of Washington, Seattle, WA 98105, United States}
\affil[10]{Rice University, Department of Chemistry, Houston, TX 77005, United States}

\begin{document}
	\title{Deeptime:~a Python library for machine learning dynamical models from time series data}
	
	\hypersetup{
		pdftitle={Deeptime: a Python library for machine learning dynamical models from time series data},
		pdfsubject={cs.LG, cs.AI, stat.ML},
		pdfauthor={Moritz Hoffmann, Martin Scherer, Tim Hempel, Andreas Mardt, Brain de Silva, Brooke E. Husic, Stefan Klus, J. Nathan Kutz, Steven L. Brunton, Frank Noé},
		pdfkeywords={python, scikit-learn, machine-learning},
	}
	
	\maketitle
	
	\begin{abstract}
		Generation and analysis of time-series data is relevant to many quantitative fields ranging from economics to fluid mechanics. In the physical sciences, structures such as metastable and coherent sets, slow relaxation processes, collective variables, dominant transition pathways or manifolds and channels of probability flow can be of great importance for understanding and characterizing the kinetic, thermodynamic and mechanistic properties of the system. \textit{Deeptime} is a general purpose Python library offering various tools to estimate dynamical models based on time-series data including conventional linear learning methods, such as Markov state models (MSMs), Hidden Markov Models and Koopman models, as well as kernel and deep learning approaches such as VAMPnets and deep MSMs. The library is largely compatible with scikit-learn, having a range of Estimator classes for these different models, but in contrast to scikit-learn also provides deep Model classes, e.g. in the case of an MSM, which provide a multitude of analysis methods to compute interesting thermodynamic, kinetic and dynamical quantities, such as free energies, relaxation times and transition paths.
		The library is designed for ease of use but also easily maintainable and extensible code. In this paper we introduce the main features and structure of the deeptime software.
	\end{abstract}
	
	\section{Introduction}
	\extrafootertext{Work performed prior to employment at Amazon.}Deeptime is an open source Python library for the analysis of time-series data; i.e., the provided methods relate to finding relationships between instantaneous data $\mathbf{x}_t$ for some $t\in[0, \infty)$ and corresponding future data $\mathbf{x}_{t+\tau}$ for some so-called \emph{lag-time} $\tau>0$. Most of the implemented methods try to estimate the behavior of processes when going from $\mathbf{x}_t$ to $\mathbf{x}_{t+\tau}$ by predicting the latter based on the former. The API is similar to what is implemented in the well-known software package scikit-learn~\cite{sklearn_api} and there is basic compatibility of the methods in the two packages via duck typing\footnote{Duck typing refers to the objects' behavior defining in which contexts it can be used, not so much its concrete type.}. Deeptime has two main goals: (1) making methods which were developed in different communities (such as molecular dynamics and fluid dynamics) accessible to a broad user base by implementing them in a general-purpose way, and (2) easy to extend and maintain due to modularity and very few hard dependencies.
	
	Deeptime offers the following main groups of methods:
	\begin{itemize}
		\item \textit{Dimension reduction of dynamical data.} One vitally important ingredient to understanding high-dimensional data is projecting them onto a low-dimensional manifold which preserves the ``interesting'' parts of the signal. One prominent linear method which can perform this task is principle component analysis (PCA)~\cite{pearson1901liii,hotelling1933analysis}. While PCA is a widely implemented method, it is not designed for extracting dynamically relevant features from time-series data. Time-lagged independent component analysis (TICA)~\cite{molgedey1994separation,naritomi2011slow} and dynamic mode decomposition (DMD)~\cite{schmid2010dynamic,Tu2014jcd} provide dimension reduction along with a best-fit linear model. Deeptime offers a range of methods which are based on the mathematical framework of transfer operators~\cite{koopman1931hamiltonian,gaspard_1998,klus2016numerical,klus2018data,Brunton2021koopman}, enabling users to study in particular kinetic properties of the data as well as find temporally metastable and coherent regions.
		
		\item \textit{Nonlinear dimension reduction.} While linear methods are widely used due to their simplicity, the  high-dimensional data manifold might not always be structured in a way in which important processes are linearly separable. For that reason, deeptime offers some featurizations, explicitly defined basis functions, and kernel methods. 
		
		\item \textit{Deep dimension reduction.} For nonlinear dimension reduction---especially in the high-data regime---there is also a weak dependency~(i.e., no strict requirement for installation) to PyTorch~\cite{NEURIPS2019_9015}, enabling the use of deep learning techniques for dimension reduction of time-series data. A variety of such methods has recently been developed, for example time-lagged Autoencoders~\cite{wehmeyer2018timelagged}, linearly-recurrent Autoencoder networks~\cite{Otto2019siads}, VAMPNets~\cite{mardt_vampnets_2018}, deep generative Markov state models~\cite{wu2018deep,mardt2021progress}, deep Koopman networks~\cite{Lusch2018natcomm}, and variational dynamical encoders~\cite{hernandez2018variational}. Some of the mentioned methods are implemented in this software.
		
		\item \textit{Markov state models (MSMs).} MSMs~\cite{schutte_direct_1999,swope_describing_2004,singhal_using_2004,noe_hierarchical_2007,noe_probability_2008,noe_constructing_2009,prinz_markov_2011,chodera_markov_2014,husic_markov_2018} are stochastic models describing temporal transitions between states in chains of events where each event only depends on its predecessor and has no dependency on events further in the past (known as the Markov property). They also fit into the mathematical framework of transfer operators. Based on MSMs one can estimate in particular kinetic properties from data.   
		
		\item \textit{Hidden markov models (HMMs).} HMMs~\cite{baum1966statistical, rabiner_tutorial_1989} are a type of model consisting of a hidden (i.e., not observable) Markov process emitting an observable output process depending on the hidden process. In comparison to MSMs, HMMs are more expressive and can produce good results where MSMs would not, but are harder to estimate.\footnote{A more effective/efficient model for hidden Markov processes with discrete output probability distributions is the observable operator model MSM~\cite{nuske_markov_2017} that can also be found within the deeptime package.}
		
		\item \textit{Sparse identification of nonlinear dynamics (SINDy).} SINDy~\cite{brunton2016sindy} identifies nonlinear governing equations with as few terms as possible from a library of candidate terms that best fit the data. In that way, it complements the dimension-reduction techniques. In particular, while most methods model and analyze the relationships of time-shifted pairs of data, SINDy predicts maps yielding the infinitesimal expected temporal change of the system's current state. On the other hand, SINDy can also predict discrete-time maps by directly relating the system's future state to the system's current state~(see Section~\ref{sec:sindy} for details).
	\end{itemize}
	
	Deeptime currently focuses on the domain-agnostic estimation of dynamical models and their analysis in terms of physically relevant quantities describing equilibrium or nonequilibrium behavior. The aim of Deeptime is not to provide tools specific for a single domain, such as molecular dynamics, but it can be easily combined with python packages that, e.g., load and featurize domain-specific data files in order to prepare such data for analysis with Deeptime~\cite{roe2013ptraj,romo2014lightweight,mcgibbon_mdtraj_2015,michaud2011mdanalysis,nguyen2016pytraj,gowers2019mdanalysis}. Alternatively there also exist time-series analysis packages that are more domain-specific~\cite{beauchamp2011MSMBuilder2,scherer_pyemma_2015,mardt_vampnets_2018,wehmeyer_introduction_2019} or implement a subset of deeptime's methods but with a wider range of options and/or more flexibility~\cite{Demo18pydmd,hmmlearn,desilva2020}. 
	As the dynamical model and its properties take the center stage in Deeptime, its aim is also \textit{not} to perform time-series forecasting, e.g., for weather or financial data, or clustering, regression, and annotation directly on the dynamical data itself. 
	For these types of tasks there is, e.g., the sktime project~\cite{loning2019sktime}\footnote{sktime also provides a curated overview of various projects dealing with time-series data: \url{https://www.sktime.org/en/latest/related_software.html}}.

	\section{Design and implementation}
	
	Deeptime is mainly implemented in and available for Python 3.7+ (as of now\footnote{\today}) and available for all major operating systems via the Python package index~(PyPI) and conda-forge~\cite{conda_forge_community_2015_4774216} . Some computationally expensive calculations are implemented in C++ using pybind11~\cite{pybind11} or if appropriate using NumPy~\cite{harris2020array} and SciPy~\cite{2020SciPy-NMeth}.
	
	The API itself is inspired (and largely compatible with) the one used by scikit-learn~\cite{sklearn_api}. In particular, deeptime offers \mintinline{python}{Estimator} classes, which can be fit on data. An important point at which deeptime's implementation is different to what is offered by scikit-learn is the following:~a call to \mintinline{python}{fit} leads to the creation of a \mintinline{python}{Model} instance; in particular, estimators can be fit multiple times and each time produce an independent model instance (therefore are model factories).  Regarding the structure of data they store, \mintinline{python}{Model}s carry the estimation results and are rather simple classes, that are akin to Python dictionaries. If possible, estimators offer a \mintinline{python}{partial_fit} method that allows the user to continuously update a model with a stream of data. This is particularly useful if the dataset does not fit into the computer's main memory. Additionally, \mintinline{python}{Model}s may also be \mintinline{python}{Transformer}s, meaning they can \mintinline{python}{transform} data based on the state of the \mintinline{python}{Model} instance. In such cases the corresponding \mintinline{python}{Estimator} also implements the \mintinline{python}{Transformer} interface, dispatching the call to the latest estimated model.
	
	In comparison, in scikit-learn an \mintinline{python}{Estimator} is also a \mintinline{python}{Model} and the estimation results are dynamically attached to the estimator instance. 
	Given that our models come with a large variety of attached methods and properties, we deviate from this paradigm to ensure clarity and component separation and to avoid an overcrowded interface.
	Furthermore, as our \mintinline{python}{Model}s are relatively lightweighted objects that are divorced from the data they have been trained on, it is straightforward to use the Python pickle module for serialization. This way, \mintinline{python}{Estimator} instances can be re-used on existing models without side effects, fostering deeptime's applicability to parameter studies. 
	
	The number of dependencies is kept as low as possible to reduce maintenance efforts. The base functionality of deeptime only depends on the established packages NumPy~\cite{harris2020array}, Scipy~\cite{2020SciPy-NMeth}, and scikit-learn~\cite{sklearn_api}.
	Dependencies to plotting routines (matplotlib~\cite{Hunter:2007}) and deep learning components (PyTorch~\cite{NEURIPS2019_9015}) are optional.
	
	The code is hosted on GitHub (\url{https://github.com/deeptime-ml/deeptime}) and licensed under LGPLv3, meaning it uses a license with weak copyleft so the library can be used also in proprietary codes. The repository is coupled to the continuous integration service Azure Pipelines, performing automated testing upon changes or proposed changes to the main branch. The project uses the pytest testing framework~\cite{pytest6.2}. 
	
	The documentation aims for maximal transparency with respect to the implemented methods and the implementation details. To that end, the main methods and their basic usage are explained in Jupyter notebooks~\cite{jupyter} with some theoretical background, references, and illustrative examples. The detailed API documentation is generated directly from the Python source code, so that it can be referred to while using the software but also while developing new components or fixing bugs. Furthermore, there are short example scripts for the datasets and selected methods, compiled into example galleries. All this is rendered into HTML and transparently hosted on GitHub pages using Sphinx under \href{https://deeptime-ml.github.io/}{https://deeptime-ml.github.io/}.
	
	The deeptime library is structured in such a way that the entire user interface is exposed at package-level. We structure the (sub-)packages as follows:
	\begin{itemize}
		\item \mintinline{python}{deeptime.base}: Contains all the basic classes of deeptime, in particular the interface definitions for \mintinline{python}{Estimator}s, \mintinline{python}{Model}s, and \mintinline{python}{Transformer}s.
		\item \mintinline{python}{deeptime.basis}: A set of basis functions which can be used for SINDy and some of the dimension reduction algorithms as ansatz and/or featurization.
		\item \mintinline{python}{deeptime.kernels}: A set of predefined kernels which can be used in kernel methods. Some of these possess subclasses with a \mintinline{python}{Torch} prefix, indicating that they are PyTorch-ready and support batched evaluation as well as backpropagation.
		\item \mintinline{python}{deeptime.sindy}: Contains an implementation of the SINDy estimator (see Section~\ref{sec:sindy}).
		\item \mintinline{python}{deeptime.covariance}: Methods for estimating covariance and autocorrelation matrices from time-series data in an online fashion. These are mainly used by some of the decomposition methods (see Section~\ref{sec:dimredux-shallow}).
		\item \mintinline{python}{deeptime.decomposition}: Decomposition methods for time-series data (see Section~\ref{sec:dimredux} for a comprehensive list of implemented estimators).
		\item \mintinline{python}{deeptime.markov}: Analysis tools, validators, and estimators for MSMs and HMMs (see Section~\ref{sec:msm}).
		\item \mintinline{python}{deeptime.clustering}: A collection of clustering/discretization algorithms. These are mostly intended for assigning frames to discrete states (potentially after using one of the dimension reduction algorithms) and subsequently estimating MSMs or HMMs.
		\item \mintinline{python}{deeptime.numeric}: A collection of numerical utilities, most notably for eigenvalue problems and regularized inverses of symmetric positive semi-definite matrices.
		\item \mintinline{python}{deeptime.data}: A selection of example data on which the algorithms can be tested (see Section~\ref{sec:data}).
		\item \mintinline{python}{deeptime.util.data}: Utilities which relate to data processing, e.g., time-series specific \mintinline{python}{DataSet} implementations which can be used in conjunction with PyTorch.
	\end{itemize}
	
	Some of the implementations are based on the molecular-dynamics analysis package PyEMMA~2~\cite{scherer_pyemma_2015,wehmeyer_introduction_2019} including its dependencies bhmm~\cite{chodera2011bayesian} and msmtools\footnote{\url{https://github.com/markovmodel/msmtools}}---modified so that they are no longer dependent on any molecular-dynamics specific libraries and offer greater flexibility---and on the dynamical systems toolbox d3s\footnote{\url{https://github.com/sklus/d3s}}. The \mintinline{python}{deeptime.sindy} package is based on and compatible to PySINDy~\cite{desilva2020}. 
	The \mintinline{python}{deeptime.decomposition} package contains an implementation of dynamic mode decomposition (DMD)~\cite{schmid2010dynamic,Tu2014jcd,Rowley2009jfm,Kutz2016book}. For a richer feature set and different variants and flavors of DMD we recommend the PyDMD package~\cite{Demo18pydmd}.
	
	\section{Dimension reduction and decomposition methods}\label{sec:dimredux}
	
	Deeptime offers a range of methods that can be used to reduce the dimension of observed data by projecting it onto dominant processes. This relates to the mathematical framework of transfer operators~\cite{koopman1931hamiltonian,gaspard_1998,mezic2005spectral,koltai2018optimal,klus2018data,wu2020variational}. We regard all operators that describe the temporal evolution of, e.g., probability densities or observables of the system's state as transfer operators. The operators we consider here are all linear operators (although in general not finite-dimensional).
	
	For an introduction to these operators we follow the presentation of~\cite{koltai2018optimal}. We distinguish two different cases: time-homogeneous processes, which possess transition probabilities that do not depend on a particular point in time (this is the case for, e.g., autonomous differential equations) and the more general case of time-inhomogeneous processes.
	
	\paragraph{Time-homogeneous processes.}
	Let $\{\mathbf{x}_t\}_{t\geq 0}$ be a Markovian and time-homogeneous stochastic process in state space $\mathbf{x}_t\in\Omega\subset\mathds{R}^d$ with transition density
	\begin{align}\label{eq:transition-density}
		p_{s,t} : \Omega\times\Omega\to\mathds{R}_{\geq 0},\quad \mathds{P}[\mathbf{x}_t\in B\mid \mathbf{x}_s = x] = \int_B p_{s,t}(\mathbf{x},\mathbf{y})\mathrm{d}\mathbf{y},
	\end{align}
	which is the probability of finding state $\mathbf{x}_t$ in a measurable set $B\subset\Omega$ given state $\mathbf{x}$ at time $s$. Time-homogeneity means that $p_{s,t}$ only depends on a lag-time $\tau := t-s$ but not on specific start and end times $s$ and $t$ individually, i.e.,
	\begin{align}\label{eq:transition-densitiy-homogeneous}
	    p_{s,t}(\mathbf{x},\mathbf{y}) = p_\tau(\mathbf{x},\mathbf{y}).
	\end{align}
	However, this does not mean that the law (or distribution) of the process
	\[ B \mapsto \law(\mathbf{x}_t)[B] := \mathds{P}[\mathbf{x}_t\in B] \]
	for sets $B\subset\Omega$ is time-independent. For example Brownian motion is a time-homogeneous process, however its law for a single particle at initial time is given by a delta peak in the initial position and converges to a uniform spatial distribution over time.
	
	Generally speaking, transfer operators describe the effect of the underlying dynamics on functions of the state $\mathbf{x}_t$.
	A particularly important transfer operator, the Koopman operator (first introduced in~\cite{koopman1931hamiltonian}), is defined as
	\begin{align} \label{eq:koopman}
		\mathcal{K}_\tau &: L^\infty(\Omega) \to L^\infty(\Omega),\quad  [\mathcal{K}_\tau g](\mathbf{x}) := \int g(\mathbf{y})p_\tau(\mathbf{x},\mathbf{y})\mathrm{d}\mathbf{y} = \mathds{E}[g(\mathbf{x}_{t+\tau})\mid \mathbf{x}_t = \mathbf{x}], 
	\end{align}
	evolving the observable $g$ for a lag-time $\tau > 0$. The function space\footnote{Strictly speaking $L^p$ consists of equivalence classes of measurable functions where the equivalence relation is defined by functions being equal ``almost everywhere'', i.e., can differ on sets of measure zero.} $g\in L^\infty(\Omega)$ is of the family of $L^p$ spaces with
	\begin{align*}
	    L^p(\Omega) := \left\{ f : \Omega\to\mathds{C} \text{ s.t. } f \text{ measurable }\land \|f\|_p := \left(\int_\Omega |f|^p \right)^{1 / p} < \infty \right\}
	\end{align*}
	for $ 1 \le p < \infty $ and $L^{\infty}(\Omega) := \{ f : \Omega\to\mathds{C} \text{ s.t. } f \text{ measurable }\land \exists C\geq 0 : |f(x)|\leq C \text{ a.e.} \}$.
	In case of deterministic dynamics $\mathbf{x}_{t+\tau} = \bm{\Psi}(\mathbf{x}_t)$, the transition density consists of delta peaks and the Koopman operator is simply the composition $\mathcal{K}_\tau g = g \circ \bm{\Psi}$.
	
	Another commonly used transfer operator to describe Markovian dynamics is the Perron--Frobenius (PF) operator~\cite{lasota1994chaos,boyarsky_laws_1997}
	\begin{align}\label{eq:perron-frobenius}
		\mathcal{P}_\tau &: L^1(\Omega) \to L^1(\Omega), \quad [\mathcal{P}_\tau f](\mathbf{y}) = \int \mathbf{f}(\mathbf{x})p_{\tau}(\mathbf{x},\mathbf{y})\mathrm{d}\mathbf{x},
	\end{align}
    which evolves probability density functions $f\in L^1(\Omega)$. Since it is a Markov operator ($\mathcal{P}_\tau f\geq 0$ and $\|\mathcal{P}_\tau f\| = \|f\|$ for all $f\geq 0$), probability density functions are mapped to probability density functions~\cite{lasota1994chaos}.
    
    The PF operator is the adjoint of the Koopman operator~\cite{lasota1994chaos,boyarsky_laws_1997}, i.e.,
	\begin{align}\label{eq:pf-koopman-adjoint-standard}
		\langle \mathcal{P}_\tau f, g\rangle = \langle f, \mathcal{K}_\tau g\rangle \quad\forall f\in L^1(\Omega), g\in L^\infty(\Omega),
	\end{align}
	where the bracket is defined as $\langle h_1, h_2 \rangle := \int_\Omega h_1(\mathbf{x}) h_2(\mathbf{x}) \mathrm{d}\mathbf{x}$. Although $L^p$ spaces with $p\neq 2$ are not Hilbert spaces, the product of two functions $h_1\in L^p(\Omega)$ and $h_2\in L^q(\Omega)$ is integrable as long as $1/p + 1/q = 1$ for $1\leq p,q\leq\infty$.
	
    For the rest of this section we assume that there exists a stationary distribution $\mu\in L^1(\Omega)$ satisfying $\mathcal{P}_\tau\mu = \mu$. If such a stationary distribution $\mu$ exists and $\mu(x) > 0$ almost everywhere, then the time-homogeneous processes $\{\mathbf{x}_t\}_{t\geq 0}$ is ergodic and the stationary distribution is unique. Vice versa, if $\{\mathbf{x}_t\}_{t\geq 0}$ is ergodic, there exists at most one stationary distribution~\cite{lasota1994chaos}.
    
    Given the stationary distribution we can define a PF operator with respect to $\mu$ (also simply called the transfer operator),
    \begin{align}\label{eq:transfer-op}
        \mathcal{T}_\tau &: L^1(\Omega) \to L^1(\Omega), \quad [\mathcal{T}_\tau u](\mathbf{y}) = \frac{1}{\mu(\mathbf{y})}\int \mu(\mathbf{x})u(\mathbf{x})p_{\tau}(\mathbf{x},\mathbf{y}) \mathrm{d}\mathbf{x}.
    \end{align}
    Instead of evolving probability densities $f$, it evolves densities $u=f/\mu$ with respect to the stationary distribution. Due to this construction we obtain the normalization $\mathcal{T}_\tau\mathds{1} = \mathds{1}$, encoding that the stationary distribution is preserved under propagation in time.
    
    Under some conditions~\cite{klus2016numerical,klus2018data,wu2020variational,tian2021kernel}, the function spaces from and to which the operators map can be assumed to be reweighted $L^2$ spaces,
	\begin{align}\label{eq:weighted-l2-space}
		L_\rho^2(\Omega) := \left\{h : \|h\|^2_\rho < \infty \;\text{with}\; \langle f, g\rangle_{\rho} := \int_\Omega f(\mathbf{x})\overline{g(\mathbf{x})}\rho(\mathbf{x})\mathrm{d}\mathbf{x} \right\},
	\end{align}
	where $\mathcal{P}_\tau : L^2_{\mu^{-1}}(\Omega) \to L^2_{\mu^{-1}}(\Omega)$, $\mathcal{T}_\tau : L^2_\mu(\Omega) \to L^2_\mu(\Omega)$, and $\mathcal{K}_\tau : L^2_{\mu}(\Omega) \to L^2_{\mu}(\Omega)$. In what follows we assume that this is the case.
	
	Via a straightforward calculation using~(\ref{eq:pf-koopman-adjoint-standard}) one obtains that Koopman operator and transfer operator are also adjoint in the reweighted spaces, i.e.,
    \begin{align}\label{eq:koopman-pf-dual}
         \langle \mathcal{T}_\tau f, g \rangle_\mu = \langle f, \mathcal{K}_\tau g\rangle_\mu \quad\forall f,g\in L^2_\mu(\Omega).
    \end{align}

	\paragraph{Time-inhomogeneous processes.} In the case of time-inhomogeneous processes, the transition density~(\ref{eq:transition-density}) depends directly on the initial and/or final time; i.e., Equation~(\ref{eq:transition-densitiy-homogeneous}) no longer holds. This also means that the operators~ (\ref{eq:koopman})--(\ref{eq:transfer-op}) no longer depend on the lag-time $\tau$ but rather on specific start and end times $s$ and $t$, respectively (equivalently: on start time $s$ and with lag-time $\tau$).
	For such systems there is in general no stationary distribution $\mu$, so we consider the distribution $\mu_s$ at initial time $s$ and $\mu_t$ at final time $t$, related by $\mu_t = \mathcal{P}_{s,t}\mu_s$.
	The transfer operator can be defined as
	\begin{align}\label{eq:transfer-op-inhomogeneous}
	    \mathcal{T}_{s,t} : L^2_{\mu_s}(\Omega)\to L^2_{\mu_t}(\Omega),\quad \mathcal{T}_{s,t} u = \frac{1}{\mu_t}\mathcal{P}_{s,t}(u\mu_s).
	\end{align}
	As in the time-homogeneous case, this operator is the adjoint of the time-inhomogeneous Koopman operator~\cite{denner2017coherent}
	\begin{align*}
	    \langle \mathcal{T}_{s,t} f, g \rangle_{\mu_t} = \langle f, \mathcal{K}_{s,t} g\rangle_{\mu_s} \quad\forall f\in L^2_{\mu_s}(\Omega)\forall g\in L^2_{\mu_t}(\Omega),
	\end{align*}
	where $\mathcal{K}_{s,t} : L^2_{\mu_t}(\Omega) \to L^2_{\mu_s}(\Omega)$.
	
	For the remainder of this section we will simplify the notation to $\mathcal{P}$, $\mathcal{T}$, and $\mathcal{K}$ for the Perron--Frobenius, transfer, and Koopman operators, respectively. Also we often wish to consider/use vector-valued feature functions, in which case it is assumed that the transfer operators act component-wise.

	One particular advantage of considering any of the transfer operators over directly analyzing the (in general highly nonlinear) temporal evolution of processes' full states is their linearity. While the considered operator usually cannot be represented as a finite-dimensional matrix, one can seek projections and/or approximations. 
	These approximations can be used to identify and project onto the slow processes as well as metastable and coherent sets~\cite{meyn2012markov,schutte2013metastability,klus2016numerical}. There are different methods available for making the approximations which vary in their assumptions, approximation power, and interpretability, some of which are accompanied by variational theorems.
	
	\subsection{Conventional dimension reduction and decomposition}\label{sec:dimredux-shallow}
	
	The conventional machine learning estimators for dimension reduction supported by deeptime are detailed below. For more thorough introductions to available methods and overviews of their relationships, we refer the reader to Refs.~\citenum{klus2018data,klus2019kernel,glielmo_unsupervised_2021}. Most of the following methods seek a matrix $K\in\mathds{R}^{m\times m}$, a finite-dimensional approximation of a transfer operator that should fulfill
	\begin{align}\label{eq:koopman-matrix}
		\mathds{E}[\mathbf{g}(\mathbf{x}_{t+\tau})] = K^\top \mathds{E}[\mathbf{f}(\mathbf{x}_t)]
	\end{align}
	as closely as possible for time series data $\mathbf{x}_t$. The system's state $\mathbf{x}_t$ is transformed into feature space by $\textbf{f},\textbf{g}\in\mathcal{F}^m$, where $\mathcal{F}$ is the space of scalar feature functions.
	
\begin{figure}
		\begin{center}
			\begin{tikzpicture}
				\node (msm) [class3]
				{
					\textbf{MSM}
				};
				\node (vac) [class3, left=2.4cm of msm, bottom color=blue!10, top color=blue!10]
				{
					\textbf{VAC}
				};
				
				\node (edmd) [class3, right=2.4cm of msm, bottom color=green!5, top color=green!5]
				{
					\textbf{EDMD}
				};
				\node (tica) [class3, above=of vac, bottom color=blue!10, top color=blue!10]
				{
					\textbf{TICA}
				};
				\node (dmd) [class3, above=of edmd, bottom color=green!5, top color=green!5]
				{
					\textbf{DMD}
				};
				\node (vamp) [class3, below=2.8cm of vac, bottom color=blue!10, top color=blue!10]
				{
					\textbf{VAMP}
				};
				\node (kcca) [class3, below=of vamp, bottom color=blue!10, top color=blue!10]
				{
					\textbf{kernel CCA}
				};
				\node (kedmd) [class3, below=of edmd, bottom color=green!5, top color=green!5]
				{
					\textbf{kernel EDMD}
				};
				\node (gmsm) [class3, right=of vamp, below=2.8cm of msm]
				{
					\textbf{GMSM}
				};
				\node (kvad) [class3, right=of gmsm, below=2.8cm of edmd, bottom color=blue!10, top color=blue!10]
				{
					\textbf{KVAD}
				};
				
				\path[pil,->] (vac) edge node[left,pos=0.5] {$\bm{\phi}(\mathbf{x})=\mathbf{x}$} (tica);
				\path[pil,->] (edmd) edge node[right,pos=0.5] {$\bm{\phi}(\mathbf{x})=\mathbf{x}$} (dmd);
				\path[pil,<->] (tica) edge node[above,pos=0.5] {dual} (dmd);
				\path[pil,<->] (vac) edge[bend left=10] node[above,pos=0.5] {equivalent} (edmd);
				\path[pil,->] (vac) edge node[below,pos=0.5] {$\phi(\mathbf{x}) = \mathds{1}_A(\mathbf{x})$} (msm);
				\path[pil,->] (edmd) edge node[below,pos=0.5] {$\phi(\mathbf{x}) = \mathds{1}_A(\mathbf{x})$} (msm);
				\path[pil,->] (vamp) edge node[left,pos=0.5] {reversible} (vac);
				\path[pil,->] (gmsm) edge node[left,pos=0.5] {} (msm);
				\path[pil,->] (vamp) edge node[left,pos=0.5] {kernelize} (kcca);
				\path[pil,->] (edmd) edge node[right,pos=0.5] {kernelize} (kedmd);
				\path[pil,->] (vamp) edge node[below,pos=0.5] {$\phi(\mathbf{x}) = \mathds{1}_A(\mathbf{x})$} (gmsm);
				\node[draw,dotted,fit=(tica) (dmd) (edmd) (vac) (kedmd) (msm),inner sep=2mm] (box1) {};
				\node[whitebox, above=0mm of box1, bottom color=white, top color=white] {\textbf{(a)} Time-homogeneous};
				\node[draw,dotted,fit=(vamp) (kvad) (kcca) (gmsm),inner sep=2mm] (box2) {};
				\node[whitebox, above=0mm of box2] {\textbf{(b)} Time-inhomogeneous};
			\end{tikzpicture}
		\end{center}
		\caption{\textbf{Relationships of conventional dimension reduction methods.} Methods shaded in green are regression-based, while methods shaded in blue are based on a variational principle. \textbf{(a)} Methods assuming time-homogeneous dynamics. While the ``variational approach to conformational dynamics''~(VAC) can only be applied under time-reversible dynamics, ``extended dynamic mode decomposition''~(EDMD) makes no assumptions about the dynamics' reversibility. \textbf{(b)} Methods supporting time-inhomogeneous dynamics. While the ``kernel embedding based variational approach for dynamical systems''~(KVAD) is based on a variational principle, the ansatz and in particular the estimated transfer operator is a Perron--Frobenius operator in contrast to the ``variational approach for Markov processes''~(VAMP), which estimates an approximation of the Koopman operator.}
		\label{fig:koopman-methods}
	\end{figure}
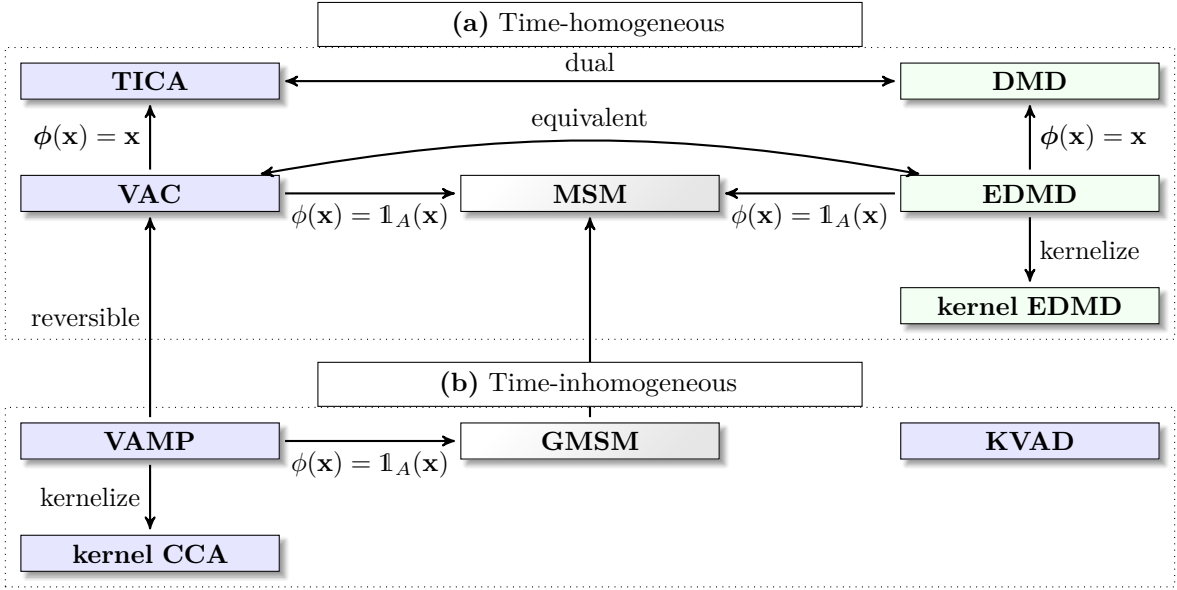

	We give an overview of conventional dimension reduction methods in Fig.~\ref{fig:koopman-methods}, all of which reside in the \mintinline{python}{deeptime.decomposition} subpackage. Roughly, the methods can be divided into groups of estimators that are restricted to data observed from time-homogeneous systems (Fig.~\ref{fig:koopman-methods}a) and estimators that are also capable of working with data of time-inhomogeneous systems (Fig.~\ref{fig:koopman-methods}b).
	
	Another distinction can be made by considering the estimation approach of the respective methods. While some are regression-based (green shade in Fig.~\ref{fig:koopman-methods}), others (purple shade) operate within the framework of an underlying variational principle.
	
	In what follows, we have instantaneous data $\mathbf{x}_i\in\mathds{R}^d$ and time-lagged data $\mathbf{y}_i\in\mathds{R}^d$ organized into matrices $X = [\mathbf{x}_1,\ldots,\mathbf{x}_n]\in\mathds{R}^{d\times n}$ and $Y = [\mathbf{y}_1,\ldots,\mathbf{y}_n]\in\mathds{R}^{d\times n}$, respectively.

	\paragraph{Dynamic mode decomposition (DMD).} DMD~\cite{schmid2010dynamic,Rowley2009jfm, Tu2014jcd,Kutz2016book} was introduced in the fluid dynamics community to extract spatiotemporal coherent structures from high-dimensional time series data. It is closely related to~(\ref{eq:koopman-matrix}) in the sense that its objective is to solve the regression problem
	\begin{align}\label{eq:dmd}
		\min_M \| Y - M_\mathrm{DMD}X \|_F
	\end{align}
	for a matrix $M_\mathrm{DMD}\in\mathds{R}^{d\times d}$. A subsequent spectral analysis of $M_\mathrm{DMD}$ can reveal information about the dominant dynamics of the system.

	There are a variety of extensions to DMD. For example, DMD algorithms have been developed that incorporate control~\cite{Proctor2016siads}, promote sparsity~\cite{Jovanovic2014pof}, are randomized~\cite{Erichson2019siads}, and act on time delay vectors~\cite{Brunton2017natcomm} (the last has a relationship to Koopman operator analysis).
	Bagheri~\cite{Bagheri2014pof} demonstrated the sensitivity of the DMD algorithm to measurement noise, motivating several noise-robust variants:~total least squares DMD~\cite{Hemati2017tcfd}, forward backward DMD~\cite{Dawson2016ef}, Bayesian DMD~\cite{Takeishi2017JCAI}, optimized DMD~\cite{Askham2018siads}, and variational DMD~\cite{azencot2019consistent}.

	While deeptime offers a basic version of DMD, most of these extensions are currently not available. The PyDMD Python package~\cite{Demo18pydmd} offers a broad range of DMD based methods.

	\paragraph{Extended dynamic mode decomposition (EDMD).}
	EDMD~\cite{williams2015data} defines a basis set of functions or observables $B := \{ \psi_1,\ldots, \psi_m\} \subset \mathcal{F}$ to construct the vector-valued function $\bm{\Psi}(x) = (\psi_1(x),\ldots,\psi_m(x))^\top\in\mathcal{F}^m$. The sought-after matrix $K$ with respect to $\textbf{f}=\textbf{g}=\bm{\Psi}$ is the solution of the regression problem
	\begin{align}\label{eq:edmd}
		\hat{K} = \argmin_K \| \bm{\Psi}(Y) - K\bm{\Psi}(X) \|_F\in\mathds{R}^{m\times m},
	\end{align}
	where the application of $\bm{\Psi}$ is column-wise.
	
	A projection onto dominant processes can be found by applying the eigenfunctions of the Koopman operator deduced from $\hat K$ and $\bm{\Psi}$ corresponding to the largest eigenvalues to the transformed input data.
	
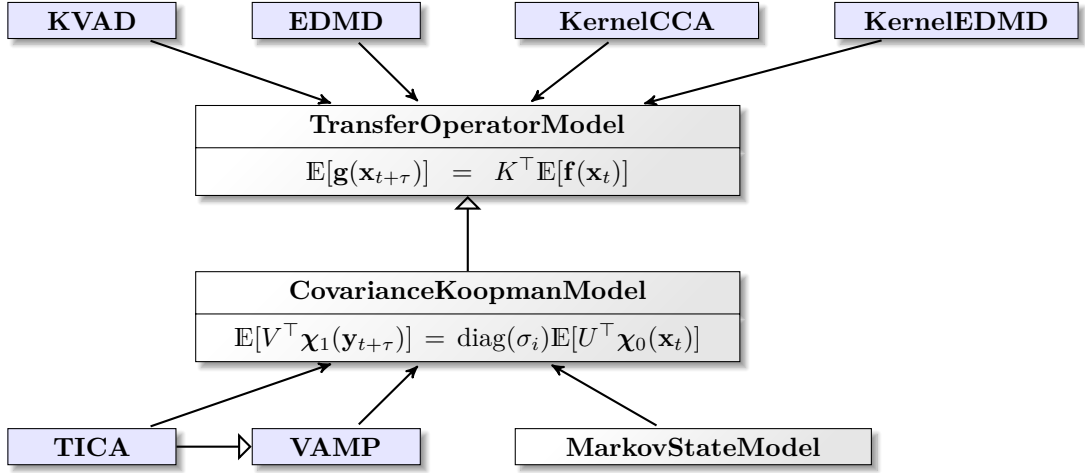
\begin{figure}
		\begin{center}
			\begin{tikzpicture}
				\node (koopmanmodel) [class, rectangle split, rectangle split parts=2]
				{
					\textbf{TransferOperatorModel}
					\nodepart{second}$\mathds{E}[\mathbf{g}(\mathbf{x}_{t+\tau})] = K^\top \mathds{E}[\mathbf{f}(\mathbf{x}_t)]$
				};
				\node (covkoopmanmodel) [class, rectangle split, rectangle split parts=2, below=of koopmanmodel]
				{
					\textbf{CovarianceKoopmanModel}
					\nodepart{second}$\mathds{E}[V^\top\bm{\chi}_1 (\mathbf{y}_{t+\tau})] = \mathrm{diag}(\sigma_i) \mathds{E}[U^\top\bm{\chi}_0(\mathbf{x}_{t})]$
				};
				\node (aux) [above=of koopmanmodel, align=center] {};
				\node (kcca) [estimator, text width=3cm, right=0.5cm of aux]
				{
					\centering\textbf{KernelCCA}
				};
				\node (edmd) [estimator, text width=2cm, left=0.5cm of aux]
				{
					\textbf{EDMD}
				};
				\node (kvad) [estimator, text width=2cm, left=of edmd]
				{
					\textbf{KVAD}
				};
				\node (kedmd) [estimator, text width=3cm, right=of kcca]
				{
					\textbf{KernelEDMD}
				};
				\node (auxlow) [below=of covkoopmanmodel] {};
				\node (vamp) [estimator, text width=2cm, left=0.5cm of auxlow]
				{
					\textbf{VAMP}
				};
				\node (tica) [estimator, text width=2cm, left=of vamp]
				{
					\textbf{TICA}
				};
				\node (msm) [class, text width=4.5cm, right=0.5cm of auxlow]
				{
					\textbf{MarkovStateModel}
				};
				\path[inharrow,->] (covkoopmanmodel) edge (koopmanmodel);
				\path[inharrow,->] (tica) edge (vamp);
				\path[pil,->] (kedmd) edge (koopmanmodel);
				\path[pil,->] (kcca) edge (koopmanmodel);
				\path[pil,->] (kvad) edge (koopmanmodel);
				\path[pil,->] (edmd) edge (koopmanmodel);
				\path[pil,->] (tica) edge (covkoopmanmodel);
				\path[pil,->] (vamp) edge (covkoopmanmodel);
				\path[pil,->] (msm) edge (covkoopmanmodel);
			\end{tikzpicture}
		\end{center}
		\caption{\textbf{Class diagram illustrating relationships between estimators producing approximations of transfer operators.} Estimators have a blue background while models are shaded in gray. The \mintinline{python}{TransferOperatorModel} implements observable transforms $f(\cdot)$ and $g(\cdot)$ as well as propagation of $f(x_t)$ with the Koopman matrix $K$. It is produced ({\protect\tikz[baseline=-0.5ex]\protect\draw[pil] (0mm,0mm) -- +(4mm,0mm);}) by \mintinline{python}{EDMD}, \mintinline{python}{KernelEDMD}, \mintinline{python}{KernelCCA}, and \mintinline{python}{KVAD} estimators. The the \mintinline{python}{CovarianceKoopmanModel} extends ({\protect\tikz[baseline=-0.5ex]\protect\draw[inharrow] (0mm,0mm) -- +(4mm,0mm);}) the \mintinline{python}{TransferOperatorModel}. It assumes the estimation to be based on covariance matrices and defines the Koopman matrix in a whitened space, where $\chi_0$ and $\chi_1$ are basis transformations of the state $x_t$ and $x_{t+\tau}$, respectively, and $U$ and $V$ are basis transform matrices. The Koopman matrix is then a diagonal matrix. The \mintinline{python}{CovarianceKoopmanModel} can be produced by \mintinline{python}{TICA}, \mintinline{python}{VAMP}, and \mintinline{python}{MarkovStateModel}s and additionally possesses a \mintinline{python}{score()} function.}
		\label{fig:koopman-inheritance}
	\end{figure}

	The solution of the regression problem~(\ref{eq:edmd}) is an approximate version of the desired property~(\ref{eq:transfer-op}) for specific choices of $\mathbf{f}$ and $\mathbf{g}$. In deeptime this is implemented by the model of the EDMD estimator being a \mintinline{python}{TransferOperatorModel} (see Fig.~\ref{fig:koopman-inheritance}).
	
	As its name suggests, DMD can be understood as a special case of EDMD in which the feature basis contains only the identity function, i.e., $\bm{\Psi}(\mathbf{x}) = \mathbf{x}$. If we define the set of basis functions to contain indicator functions for a given discretization of the state space, EDMD estimates MSMs (see Fig.~\ref{fig:koopman-methods} and Section~\ref{sec:msm} for details on MSMs).
	
	\paragraph{Time-lagged independent component analysis (TICA).} TICA~\cite{molgedey1994separation} is a linear transformation method which was introduced for molecular dynamics in~\cite{naritomi2011slow}, was independently derived as a method for extracting the slow molecular order parameters by invoking the variational approach for conformation dynamics (see below)~\cite{perez-hernandez_identification_2013}, and introduced as a method for constructing high-accuracy MSMs in~\cite{perez-hernandez_identification_2013,schwantes2013improvements}.
	TICA is designed for time-homogeneous processes and also assumes that the process is reversible, although it may still perform well practically when applied outside these constraints. A process is defined to be reversible if it fulfills the detailed balance condition
	\begin{align}\label{eq:detailed-balance}
	    \mu(\mathbf{x})p_\tau(\mathbf{x},\mathbf{y}) = \mu(\mathbf{y})p_\tau(\mathbf{y},\mathbf{x})\quad\forall \mathbf{x},\mathbf{y}\in\Omega.
	\end{align}
	As a consequence, the transfer~(\ref{eq:transfer-op}) and Koopman~(\ref{eq:koopman}) operators are identical and therefore self-adjoint. Assuming $\mathcal{K}$ to be a Hilbert--Schmidt operator, this means that (using the Hilbert--Schmidt theorem) there is an eigenvalue decomposition
	\begin{align}
	    \mathcal{K} = \sum_{i=1}^\infty\lambda_i\langle \cdot, \varphi_i\rangle_\mu\varphi_i,
	\end{align}
	where $\varphi_i$ are eigenfunctions with $\langle\varphi_i,\varphi_{i'}\rangle_\mu = \delta_{ii'}$ and $\lambda_i$ are eigenvalues which are real and bounded by the eigenvalue $\max_i\lambda_i=1$ with a multiplicity of one (see~\cite{nuske2014variational}).

	The objective of TICA is to yield components which are uncorrelated and also maximize the time-autocorrelation under lag-time $\tau$. To this end, one can solve the generalized eigenvalue problem
	\begin{align}\label{eq:tica-eigenvalue-problem}
	    C_{0\tau}\hat\varphi_i = \hat\lambda_iC_{00}\hat\varphi_i,
	\end{align}
	where $C_{00} = \frac{1}{n-1}XX^\top$ is the instantaneous covariance matrix and $C_{0\tau} = \frac{1}{n-1}XY^\top$ is the time-lagged covariance matrix. The reversibility assumption leads to $C_{00} = C_{\tau\tau}$ and $C_{0\tau}=C_{\tau 0} = C_{0\tau}^\top$ and therefore eigenvalues $\hat{\lambda}_i\in\mathds{R}$. Because real numbers possess a total order, we can assume that the eigenvalues $\hat\lambda_i$ are in a descending order and the transformation $\hat{\bm{\varphi}}(\cdot) = [\hat{\varphi}_1(\cdot),\ldots,\hat{\varphi}_k(\cdot)]$ is the TICA projection onto the first $k$ dominant components. The corresponding eigenvalues can be related to relaxation timescales of the processes $\hat{\varphi}_i$~\cite{perez-hernandez_identification_2013}. Therefore, if we know a priori that the system is time-homogeneous and reversible, TICA can be more data-efficient and yield more interpretable results compared to methods, which do not make these assumptions.
	
	For a comparison with DMD it is useful to identify TICA with $M_\mathrm{TICA} = C_{00}^\dagger C_{0\tau}$, where $C_{00}^\dagger$ denotes the Moore--Penrose pseudoinverse. It can be seen that $M_\mathrm{TICA} = M_\mathrm{DMD}^\top$~\cite{klus2018data}. Therefore, the TICA transformation consists of Koopman eigenfunctions projected onto the basis spanned by $\psi_k(\mathbf{x})=\mathbf{x}_k$ and the DMD modes are the corresponding Koopman modes---, i.e., the coefficients $\eta_k = (\eta_{k1},\ldots,\eta_{kd})^\top$ required to write the $k$-th component of the full-state observable in terms of eigenfunctions $g_k(\mathbf{x}) = \mathbf{x}_k = \sum_i\eta_{ki}\varphi_i(\mathbf{x})$~\cite{williams2015data,klus2018data}.
	
	This relationship is reflected in Fig.~\ref{fig:koopman-methods} by identifying DMD and TICA as ``dual''. This duality can also be found within the deeptime software:~in contrast to DMD, TICA is a subclass of \mintinline{python}{TransferOperatorModel}~(see Fig.~\ref{fig:koopman-inheritance}).
	
	\paragraph{Variational approach for conformational dynamics (VAC).} Like TICA,  VAC~\cite{noe2013variational,nuske2014variational} assumes time-homogeneous and reversible dynamics. Similar to the generalization from DMD to EDMD, VAC generalizes TICA using a basis $B := \{\psi_1, \ldots, \psi_m\}\subset\mathcal{F}$ to construct a transformation $\bm\Psi(\mathbf{x}) = (\psi_1(\mathbf{x}),\dots,\psi_m(\mathbf{x}))^\top$. Subsequently the instantaneous and time-lagged data is transformed to $\bm{\Psi}(X)$ and $\bm{\Psi}(Y)$, respectively, and used in the TICA problem instead of $X$ and $Y$. From this it becomes clear that TICA can be understood as a special case of VAC with $\bm{\Psi}(\mathbf{x}) = \mathbf{x}$~(see Fig.~\ref{fig:koopman-methods}). Because it is algorithmically identical to TICA under a prior featurization of data, there is no dedicated VAC estimator in deeptime. Under the particular choice of basis functions being indicator functions, VAC estimates MSMs~(see Fig.~\ref{fig:koopman-methods} and Section~\ref{sec:msm} for details on MSMs).
	
	As its name suggests, VAC involves a variational bound. It defines the score $s_\mathrm{VAC} := \sum_i\hat{\lambda}_i$ which is bounded from above by the sum over the eigenvalues of the true Koopman operator and therefore expresses how much of the slow dynamics is captured in the projection~\cite{noe2013variational,nuske2014variational}. The score can be used to optimize the feature functions $\bm{\Psi}$. We will see in the following paragraph that under the assumption of reversible dynamics, the VAC score is equal to the VAMP-$1$ score, which is why deeptime only offers a VAMP score implementation. Assuming reversible dynamics, VAC is equivalent to EDMD~(see Fig.~\ref{fig:koopman-methods}).
	
	\paragraph{Variational approach for Markov processes (VAMP).} VAMP~\cite{wu2020variational}, sometimes also referred to as ``time-lagged canonical correlation analysis'' (TCCA)~\cite{husic2019deflation}, not only optimizes for $K$ but also optimizes for $\mathbf{f}$ and $\mathbf{g}$. This cannot be achieved by merely solving the regression problem~(\ref{eq:edmd})---as, e.g., the trivial model $\mathbf{f} = \mathbf{g} \equiv (1,\ldots,1)^\top$, $K=\mathrm{Id}$ is not informative but yields zero error.
	Instead, VAMP minimizes the left-hand side of
	\begin{align}\label{eq:koopman-approx-error}
		\|\mathcal{K} - \hat{\mathcal K}\|_{\mathrm{HS}}^2 = -\mathcal{R}(\mathbf{f}, \mathbf{g}) + \|\mathcal{K}\|_{\mathrm{HS}}^2,
	\end{align}
	the Hilbert--Schmidt norm of the difference between true Koopman operator~(\ref{eq:koopman}) and approximated Koopman operator $\hat{\mathcal K}$ deduced from $K$, $\mathbf{f}$, and $\mathbf{g}$. The minimization is achieved by maximizing $\mathcal{R}$, a variational score. The decomposition~(\ref{eq:koopman-approx-error}) of the modelling error assumes that $\mathcal{K}$ is indeed a Hilbert--Schmidt operator.
	
	In~\cite{wu2020variational} it was shown that the smallest approximation error~(\ref{eq:koopman-approx-error}) is achieved for
	\begin{align}
		\hat{\mathcal{K}} = \sum_{i=1}^m\sigma_i\langle\cdot, \phi_i\rangle\psi_i,
	\end{align}
	where $K = \mathrm{diag}(\sigma_1,\ldots,\sigma_m)$, $\textbf{f}=(\psi_1,\ldots,\psi_m)$, $\textbf{g}=(\phi_1,\ldots,\phi_m)$, and  $\sigma_i,\psi_i,\phi_i$ are the square root of the $i$-th eigenvalue, left eigenfunction, and right eigenfunction of the forward-backward operator $\mathcal{K}^*\mathcal{K}$, respectively.
	
	During estimation~(similar to TICA and VAC), covariance matrices are estimated and under regularization inverted to perform whitening operations to finally construct an approximation of the Koopman operator. One obtains coefficient matrices $U,V\in\mathds{R}^{m\times k}$ and the matrix $K\in\mathds{R}^{k\times k}$, so that
	\begin{align}
		\mathds{E}[V^\top\bm{\chi}_1(\mathbf{x}_{t+\tau})] \approx K^\top \mathds{E}[U^\top\bm{\chi}_0(\mathbf{x}_{t})],
	\end{align}
	where $\mathbf{\chi}_0$ and $\mathbf{\chi}_1$ are vectors of basis functions which optimally should contain $\psi_i$ and $\phi_i$ in their span, respectively.

	The family of VAMP-$r$ scores, 
	\begin{align}\label{eq:vamp-r-score}
		\mathcal{R}_r := \sum_i\sigma_i^r,
	\end{align}
	as well as the VAMP-E score~(see Ref.~\citenum{wu2020variational} for a definition) can be optimized to minimize the model error on the left-hand side of~(\ref{eq:koopman-approx-error}) and therefore can be used to select optimal features and/or observables by using cross-validation techniques~(see, e.g.,~\cite{scherer2019variational}). These scores give rise to the ``variational'' aspect of VAMP as they are bounded from above and their maximization leads to better approximations.
	
	VAMP therefore generalizes VAC to a time-inhomogeneous and nonreversible setting~(recall Fig.~\ref{fig:koopman-methods}). While VAMP is applicable in more situations, i.e., because it possesses greater generality and nonequilibrium dynamics are more common in nature, it also loses some of its interpretability---as, e.g., singular values can in general no longer by related to relaxation timescales of processes.

	The deeptime library reflects the mathematics of the VAMP approach by the VAMP estimator producing a \mintinline{python}{CovarianceKoopmanModel}, an extension of the \mintinline{python}{TransferOperatorModel}, which in particular allows the evaluation of VAMP scores~(see Fig.~\ref{fig:koopman-inheritance}). The estimator can deal with large amounts of data, because the estimation procedure is based on the decomposition of covariance matrices, which can be constructed incrementally~\cite{chan1982updating}. Furthermore, TICA is a subclass of VAMP, as the two methods are algorithmically closely related.\footnote{While TICA and VAC are special cases of VAMP, in deeptime
	the estimators are not combined into one due to differences in how covariance matrices are estimated---in particular, TICA's stronger inductive bias is implemented by forced symmetrizations which are not applicable to VAMP---and differences in the decomposition~(eigenvalue decomposition and singular value decomposition for TICA and VAMP, respectively).}
	
	Analogously to VAC, the choice of indicator feature functions leads to generalized MSMs (GMSMs)~\cite{koltai2018optimal}, which are also applicable to time-inhomogeneous systems~(see Fig.~\ref{fig:koopman-methods}).
	
	\paragraph{Kernel canonical correlation analysis (Kernel CCA).} Kernel CCA~\cite{bach2002kernel} is a kernelized version of canonical correlation analysis~(CCA)~\cite{hotelling1936relations} that seeks to maximize the correlation between two multidimensional random variables $X$ and $Y$~(pairs of instantaneous and time-lagged data, respectively). In kernel CCA, the standard inner products are replaced by a kernel function $\kappa(\cdot, \cdot)$ using the ``kernel trick''. Deeptime has a subpackage dedicated to kernel implementations (\mintinline{python}{deeptime.kernels}), containing (amongst others) vectorized versions of the popular Gaussian kernel 
	\begin{align}\label{eq:gaussian-kernel}
		\kappa(\mathbf{x}, \mathbf{x}') = \exp\left(-\frac{1}{2}\|\mathbf{x}-\mathbf{x}'\|_2^2/\sigma^2\right)
	\end{align}
	as well as the polynomial kernel $\kappa(\mathbf{x}, \mathbf{x}') = (c + \mathbf{x}^\top \mathbf{x}')^p$.
	
	It was shown in~\cite{klus2019kernel} that kernel CCA can be derived from optimizing the VAMP-1 score~(\ref{eq:vamp-r-score}) within a kernel approach and thus can be understood as a kernelized version of VAMP (for this reason it is sometimes referred to as kernel VAMP).
	
	In addition to the kernel parameters, the estimator also possesses a regularization parameter $\varepsilon$, as kernel CCA involves inverting covariance operators~(which on their own are generally not invertible).

	\paragraph{Kernel extended dynamic mode decomposition (Kernel EDMD).} Kernel EDMD~\cite{williams2016kernel,klus2020eigendecompositions} is, analogously to kernel CCA, a kernelized version of EDMD. In contrast to kernel CCA, it assumes a time-homogeneous process. Furthermore, kernel EDMD requires a regularziation parameter $\varepsilon$ in order to ensure invertibility of covariance operators.
	
	\paragraph{Kernel embedding based variational approach for dynamical systems (KVAD).} KVAD~\cite{tian2021kernel} is an alternative to VAMP which can also be applied to systems in which the transfer operator $\mathcal{T}$ is not Hilbert--Schmidt as an operator from $L_{\mu_s}^2$ to $L_{\mu_t}^2$\footnote{In case of time-homogeneous processes, we have $\mu_s=\mu_t=\mu$, which is the stationary distribution.}, which is~(e.g.) the case for some deterministic systems. To this end, the similarity of functions of interest is not determined using norms of $L^2$ function spaces but rather using kernel embeddings of said functions. In particular, for a given kernel $\kappa(\mathbf{x}, \mathbf{x}') = \langle \varphi(\mathbf{x}), \varphi(\mathbf{x}')\rangle$, functions $q$ can be embedded via
	\begin{align}
	\mathcal{E}q = \int\varphi(\mathbf{x})q(\mathbf{x})\mathrm{d}\mathbf{x}.
	\end{align}
	The similarity between functions $q$ and $q'$ can then be measured as
	\begin{align}\label{eq:kvad-kernel-similarity}
		\| q - q' \|_\mathcal{E} = \langle\mathcal{E} (q - q'), \mathcal{E}(q - q') \rangle.
	\end{align}
	In Ref.~\citenum{tian2021kernel} it was shown that for universal and bounded kernels $\kappa$, the Hilbert--Schmidt assumption is always fulfilled if the PF operator is considered as 
	\begin{align}\label{eq:perron-frobenius-kvad}
		\mathcal{P}_\tau : L^2_{\mu_s^{-1}} \to L^2_\mathcal{E},
	\end{align}
	where $L^2_{\mathcal{E}} = \{ f\in L^2 : \| f\|_\mathcal{E} < \infty \}$ is an $L^2$ space equipped with the kernel similarity measure~(\ref{eq:kvad-kernel-similarity}). Note that in this case the Perron--Frobenius operator as defined in~(\ref{eq:perron-frobenius-kvad}) is in general no longer the adjoint of the Koopman operator.

	Like VAMP, KVAD is based on the optimization of a~(variational) score that is bounded from above and expresses the quality of the found approximation. A key difference is the ansatz:~While VAMP yields approximations of the Koopman operator, KVAD estimates its adjoint\footnote{Adjoint in the sense of Equation~(\ref{eq:pf-koopman-adjoint-standard}).}, the PF operator. To this end, KVAD uses the transition density~(\ref{eq:transition-density}) and assumes that it can be represented as
	\begin{align}\label{eq:kvad-transition-density}
		\hat{p}_\tau(\mathbf{x}_t, \mathbf{x}_{t+\tau}) = \mathbf{f}(\mathbf{x}_t)^\top \mathbf{q}(\mathbf{x}_{t+\tau}),
	\end{align}
	where $\mathbf{q} = (q_1,\ldots,q_m)^\top$ are $m$ density basis functions and $\mathbf{f}$ are, as in~(\ref{eq:koopman-matrix}), feature functions of the system's state. This leads to the linear model~(\ref{eq:koopman-matrix}) with $\mathbf{f} = \mathbf{g}$ and
	\begin{align*}
		K = \int \mathbf{q}(\mathbf{y})\mathbf{f}(\mathbf{y})^\top \mathrm{d}\mathbf{y}.
	\end{align*}
	It has been shown~\cite{tian2021kernel} that $\mathbf{q}$ can be estimated directly from data in a nonparametric fashion, which means that all the model's parameters reside inside the definition of $\mathbf{f}$. With the help of estimated $\mathbf{f}$ and $\mathbf{q}$, also the transition matrix $K$ can be constructed. This kind of ansatz---sans the modified codomain in~(\ref{eq:perron-frobenius-kvad})---is similar to what was used in Ref.~\citenum{wu2018deep}.
	
	\subsection{Deep dimension reduction and decomposition}\label{sec:dimredux-deep}
	In addition to the conventional learning methods introduced in Section~\ref{sec:dimredux-shallow}, deeptime also offers several deep learning methods for dimension reduction.
	
	The deep learning components require PyTorch; however, PyTorch-dependent parts of the library are separate---i.e., a working installation of PyTorch is not required for the rest of the library. The estimators providing deep dimension reduction can be found in the \mintinline{python}{deeptime.decomposition.deep} subpackage.
	
\begin{figure}
		\centering
		\includegraphics[width=.6\columnwidth]{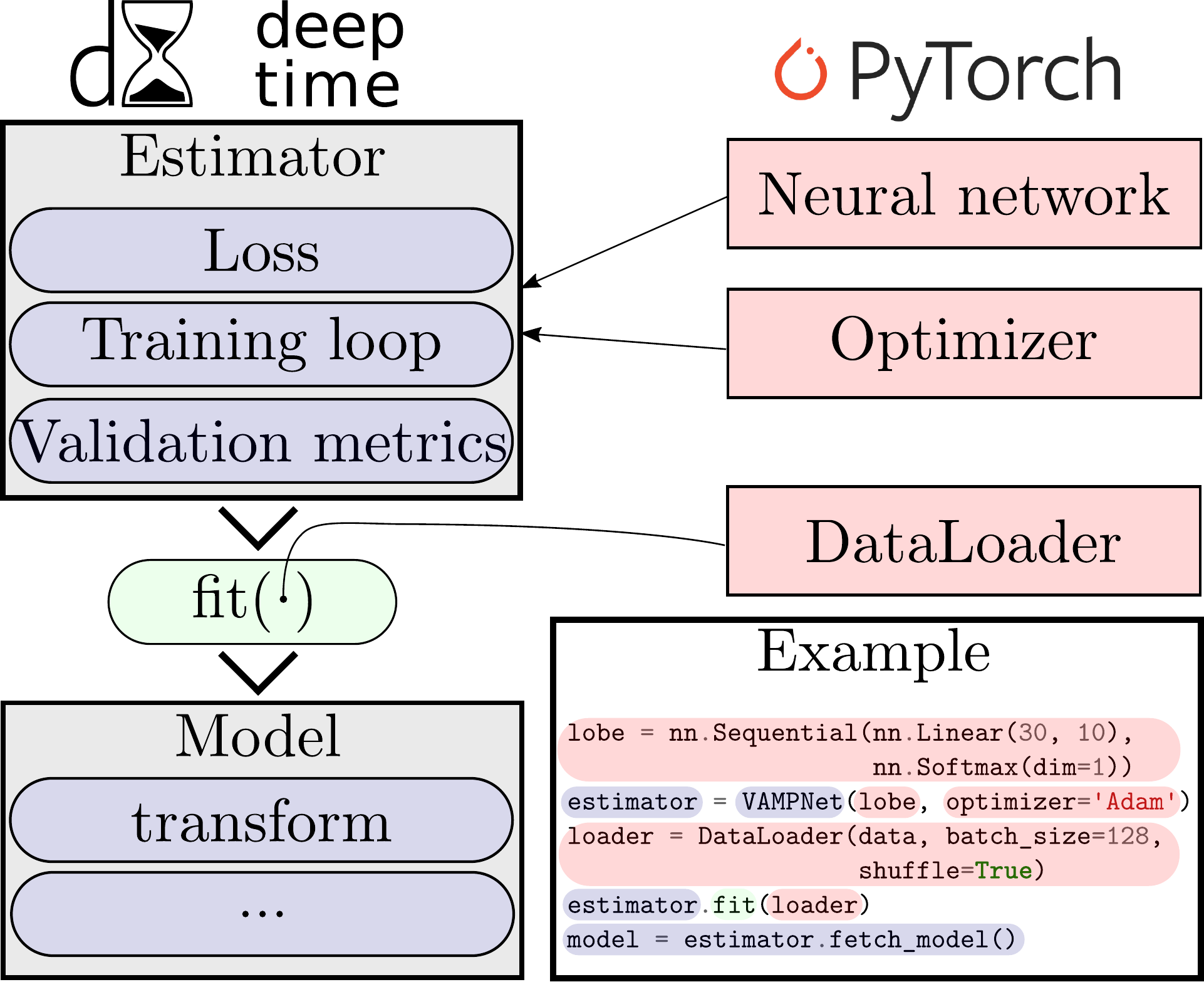}
		\caption{\textbf{Schematic overview of the flexible integration of deep learning components in deeptime.} \textbf{(a)} Deeptime interface: Some \mintinline{python}{Estimator} instances in deeptime, for example VAMPNets and Time-lagged Autoencoders, take neural networks as configurational input. These estimators offer a \mintinline{python}{fit( )} method which takes PyTorch data loaders. The data loaders are configured by the users and manage, e.g., shuffling and batch sizes. After a \mintinline{python}{Model} has been fit, it can be used for further analysis of the input data. Models containing deep learning components are designed so that they can also seamlessly work with NumPy arrays and interact with the rest of the deeptime library. \textbf{(b)} PyTorch interface: Neural networks can be defined using PyTorch and also be trained using optimizers which either are already included in PyTorch or at the very least are PyTorch compatible. Instances thereof can be used with certain estimators in deeptime. \textbf{(c)} Example: A typical workflow for a deep learning estimator, color-coded according to which library the classes and respective instances belong.}
		\label{fig:deep-overview}
	\end{figure}

	Deep learning requires some additional flexibility and data handling compared to conventional learning. In particular, the user first must define a neural network architecture and an optimizer for adjusting the network's weights. There are some predefined architectures directly available in deeptime (such as multilayer perceptrons), but in principle these as well as the optimizer are defined with PyTorch~(see Fig.~\ref{fig:deep-overview}b). Once defined, one can construct a deep estimator that contains losses, validation metrics, and training procedures~(see Fig.~\ref{fig:deep-overview}a). Fitting deep learning components typically involves shuffling and dividing the data into batches. Since the optimal batch size and also the shuffling method are problem-dependent, these choices must be made by the user. PyTorch offers \mintinline{python}{DataLoader}s for this exact purpose. Therefore \mintinline{python}{estimator.fit} is performed on a data loader instance rather than arrays~(see Fig.~\ref{fig:deep-overview}a). Finally, deep learning estimators also produce models which encapsulate among other things a copy of the trained neural network. While PyTorch neural networks operate on \mintinline{python}{torch.Tensor} instances, deeptime models with deep learning components are designed so that they can also work with ordinary NumPy arrays, ensuring a seamless integration with other and in particular conventional models and methods. For more details about PyTorch, see the official documentation\footnote{\url{https://pytorch.org/}}.
	
	\paragraph{VAMPNets.}
	VAMPNets~\cite{mardt_vampnets_2018} are a deep learning approach that seek to find parametrizations of neural networks $\bm{\chi}_0$ and $\bm{\chi}_1$ (referred to as ``lobes'') so that the VAMP-E or one of the VAMP-$r$ scores~(\ref{eq:vamp-r-score}) under these transformations is maximized, leading to smaller model errors~(recall paragraph about VAMP in Section~\ref{sec:dimredux-shallow}). This is possible because there is a variational upper bound to the scores and their computation is differentiable---therefore any of the VAMP scores can be used directly as an objective function in a deep learning context.

	Other deep learning methods which are not currently included in deeptime but also approximate the Koopman operator are, e.g., those found in Refs.~\citenum{Takeishi2017nips,Lusch2018natcomm,yeung2019learning,Otto2019siads}.
	
	\paragraph{KVADNets.}
	Analogously to VAMPNets, KVADNets optimize the variational KVAD score to find an optimal parametrization of feature functions $\mathbf{f}$~(\ref{eq:kvad-transition-density}). As in the case of VAMPNets, the KVAD score is differentiable~\cite{tian2021kernel}.
	
	\paragraph{Time-lagged (variational) autoencoders (T(V)AEs).} TAEs~\cite{wehmeyer2018timelagged} are a type of neural network approach in which instantaneous data $\mathbf{x}_t\in\mathds{R}^d$ is compressed / encoded through a parameterized function
	\begin{align*}
	    E : \mathds{R}^d \to\mathds{R}^n, \mathbf{x}_t\mapsto E(\mathbf{x}_t) = \mathbf{z}_t
	\end{align*}
	with $n \leq d$ and then reconstructed as time-lagged data $\mathbf{x}_{t+\tau}$, $\tau > 0$ via a decoder network
	\begin{align*}
	    D : \mathds{R}^n \to\mathds{R}^d, \mathbf{z}_t\mapsto D(\mathbf{z}_t) \approx \mathbf{x}_{t+\tau}.
	\end{align*}
	The optimization target is to reduce the mean-squared error between $\mathbf{x}_{t+\tau}$ and $(D\circ E)(\mathbf{x}_t)$, effectively training a latent and lower-dimensional representation $E(\mathbf{x}_t)$ of the process. In~\cite{wehmeyer2018timelagged} it was shown that in the linear case TAEs perform time-lagged canonical correlation analysis, cf. the paragraph about VAMP in Section~\ref{sec:dimredux-shallow}. An architecture that is akin to the one of TAEs was used in~\cite{chen2018molecular} to find collective variables in the context of molecular enhanced sampling.
	
	A natural extension to TAEs is to exchange the neural network architecture of an autoencoder by the architecture of a variational autoencoder (VAE)~\cite{kingma2013auto,rezende2014stochastic}, yielding the generative TVAE that can also be found in the deeptime library. In~\cite{hernandez2018variational} these architectures~(there called ``variational dynamics encoder'' (VDE)) were used in conjunction with a loss term inspired by saliency maps~\cite{kadir2001saliency}~(known from computer vision) to produce interpretable dynamical models while still maintaining the high degree of nonlinearity that can be achieved by neural networks.
	
	\subsection{Numerical experiments}
	We compare some of the dimension reduction methods introduced in Section~\ref{sec:dimredux-shallow} and Section~\ref{sec:dimredux-deep}. The first example highlights differences in the approximation if used for dimension reduction in a time-homogeneous system. The second example uses data obtained from a time-inhomogeneous system with the objective to find coherent structures.
	
	\subsubsection{Dimension reduction}
\begin{figure}
		\centering
		\includegraphics{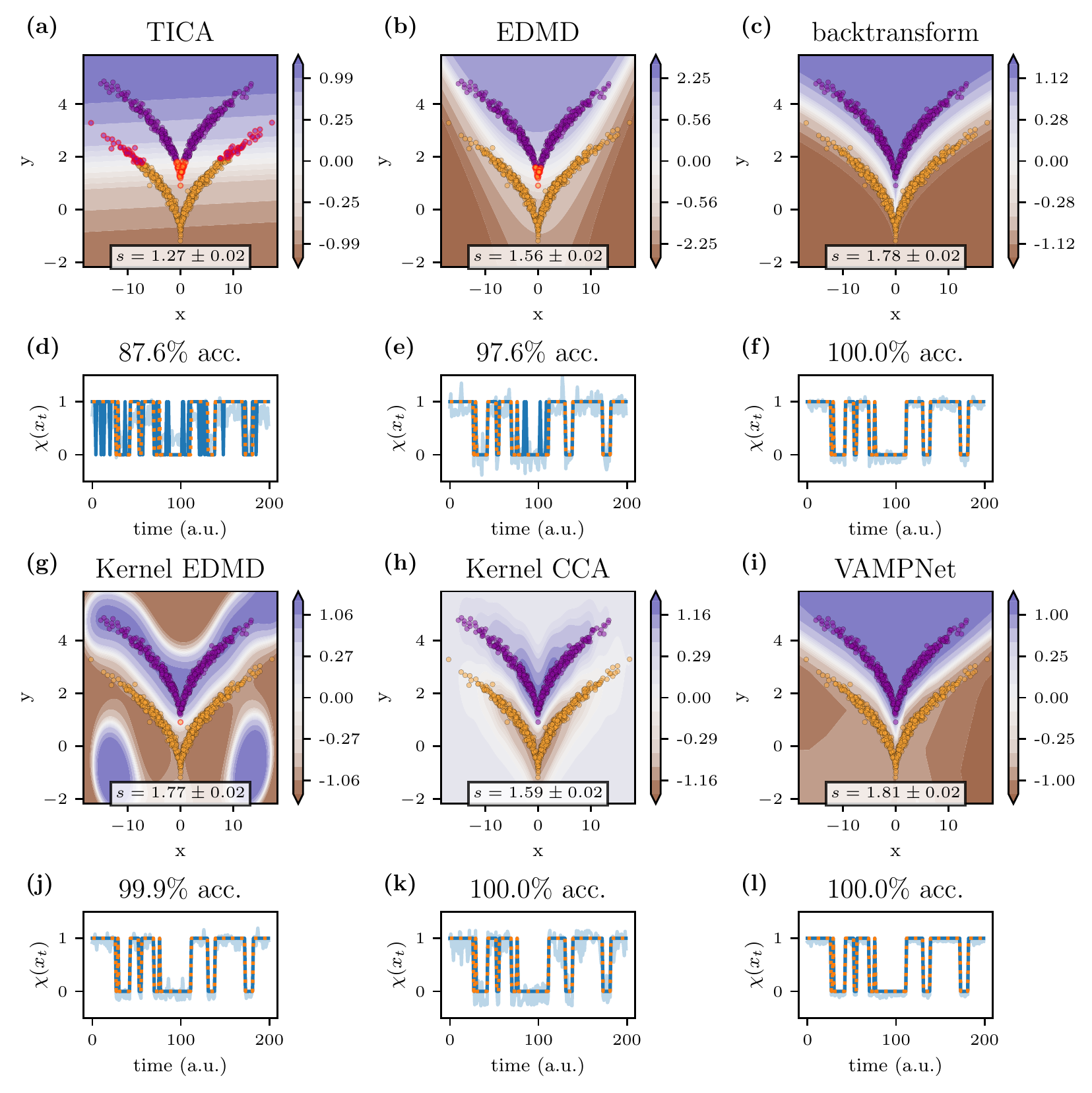}
		\caption{\textbf{Comparison of different Koopman operator methods.} The data is generated by a hidden Markov model with two hidden states and respective two-dimensional output probability distributions which are not linearly separable. The methods should approximately recover the underlying hidden process, where in \textbf{(a,b,c,g,h,i)} the background contour is the decision landscape and the scatter colors denote sharp assignments obtained by a two-state clustering in the projection. Incorrectly assigned states have a red edge. The scores reported in the bottom parts of the plots are cross-validated VAMP-2 scores. Plots \textbf{(d,e,f,j,k,l)} show a projection of the two-dimensional time-series onto one dimension (transparent blue) with crisp assignments (opaque blue) and the ground truth (orange) as reference. The accuracy~(acc.) refers to the amount of correctly assigned states after clustering.}
		\label{fig:sqrt-vamp}
	\end{figure}

	We consider a two-state hidden Markov model~(see Section~\ref{sec:msm}) with transition matrix
	\begin{align}\label{eq:sqrt-transition-matrix}
		P = \begin{pmatrix}
			0.95 & 0.05 \\ 0.05 & 0.95
		\end{pmatrix}
	\end{align}
	and anisotropic but linearly separable two-dimensional Gaussian emission distributions. In a subsequent step the data is transformed via
	\begin{align}\label{eq:sqrt-model}
		(x,y)\mapsto (x, y+\sqrt{|x|}).
	\end{align}
	This leads to two wedge-shaped output distributions which are no longer linearly separable~(see Fig.~\ref{fig:sqrt-vamp}). We simulate a trajectory of $T=1000$ frames from this model and try to recover a separation into the two original hidden states using different dimension reduction methods by projecting onto the dominant slow process,  which is the jump process between the two wedges.
	
	Fig.~\ref{fig:sqrt-vamp}~(a,b,c,g,h,i) shows the sampled and transformed time-series data $\mathbf{x}_t\in\mathds{R}^2$ as a scatter plot. The estimated models yield two-dimensional decision landscapes $\chi : \mathds{R}^2\to\mathds{R}$ which are shown in the background as a filled contour. Based on the decision landscape we obtain crisp assignments to one of the two states via $k$-means~\cite{lloyd1982least} clustering with $k=2$ cluster centers in the projected space; these clusters determine the point colors in the scatter plot. In Fig.~\ref{fig:sqrt-vamp}~(a,b,c,g,h,i) we furthermore report the $10$-fold cross-validated VAMP-$2$ score along with its standard deviation (see Section~\ref{sec:dimredux-shallow} for the score and Ref.~\citenum{mcgibbon2015variational} for the cross-validation scheme), which enables a quantitative assessment of the quality of the projection $\chi$.
	
	Fig.~\ref{fig:sqrt-vamp}~(d,e,f,j,k,l) shows the projection of the two-dimensional time series onto one dimension for the first $200$ frames of the trajectory, where $\chi(x_t)$ is presented in transparent blue with corresponding crisp clustered assignments in opaque blue and the hidden reference state is presented in orange. We also report the assignment accuracy of the crisp state with respect to the hidden reference state over the entire dataset in the titles of subfigures~(a-f).2. While the assignment accuracy can be used as another measure of the quality of the projection, it is only available if the ground truth is known. A VAMP score, on the other hand, can always be evaluated. Here, we chose the VAMP-$2$ score, because its maximization can be identified with the maximization of kinetic variance~(see Ref.~\citenum{noe2015kinetic}). Maximizing kinetic variance achieves an optimal separation of metastable sets, which corresponds to the separation of the two wedges in this example.
	
	In the limit of infinite data (i.e., when faithfully representing the original data distribution) and optimal featurization, the score should approach $s_\text{lim}=1.81$. This can be found from using the ground truth hidden transition matrix~(\ref{eq:sqrt-transition-matrix}) and applying it to the VAMP score assuming that the distribution of data is given by the stationary distribution.\footnote{In more detail, if $\bm{\mu} = (0.5, 0.5)^\top$ is the stationary distribution corresponding to~(\ref{eq:sqrt-transition-matrix}), we assume that data distribution is given by the stationary distribution and set covariance matrices $C_{00} = C_{\tau\tau} = \diag(\bm{\mu})$ and the cross-covariance matrix to $C_{0\tau} = P$~(in this example $\tau=1$). From the covariance matrices we can obtain the Koopman matrix~(cf. Section~\ref{sec:dimredux-shallow}) which can be decomposed and used for scoring.}

	Below, we discuss and further describe each of the panels in Fig.~\ref{fig:sqrt-vamp}.
	\begin{enumerate}
		\item[(a)] \textit{TICA}. TICA is a linear method in that it can only draw linear decision boundaries and the dataset is deliberately not linearly separable. Therefore the tip of the upper wedge and the outer areas of the lower wedge are misclassified. This is also reflected in the comparably low VAMP-2 score and accuracy.
		\item[(b)] \textit{EDMD}. We choose EDMD with an ansatz basis of monomials up to degree two in two-dimensional space; i.e.,
		\begin{align*}
			B = \left\{ (x,y) \mapsto x^py^q : p,q\in\mathds{N}_{\geq 0},\; p+q \leq 2 \right\}.
		\end{align*}
		This leads to a decision landscape shaped like a rounded cone, able to separate most of the data into the two hidden states except for the tip of the upper wedge. Consequently, score and accuracy achieve a higher value than the one obtained from TICA.
		\item[(c)] \textit{Backtransform}. Here we use the hand-tailored transformation $(x,y)\mapsto (x, y - \sqrt{|x|})$, which makes the two states linearly separable again and apply VAMP. This featurization uses the ground truth as prior knowledge and therefore achieves perfect state separation. Consequently, the accuracy is at $100\%$ and the VAMP-2 score reaches a high value. Due to finite data it does not quite reach the theoretical limit of $s_\mathrm{lim}=1.81$.
		\item[(d)] \textit{Kernel EDMD}. We use kernel EDMD with a Gaussian kernel~(\ref{eq:gaussian-kernel}). The regularization parameter $\varepsilon$ of the estimator as well as the bandwidth $\sigma$ of the kernel are tuned to maximize the VAMP-$2$ score on a validation set using the SLSQP optimizer~\cite{kraft1988software}, yielding $\sigma\approx 1.42$ and $\varepsilon\approx 6.7\times 10^{-4}$. The method finds a good separation between the two hidden states.
		\item[(e)] \textit{Kernel CCA}. As in the kernel EDMD case, we choose a Gaussian kernel~(\ref{eq:gaussian-kernel}) with regularization parameter and bandwidth tuned to maximize the VAMP-$2$ score on a validation set using the SLSQP optimizer~\cite{kraft1988software}. This leads to $\sigma\approx 0.85$ and $\varepsilon\approx 0.36$. Compared to the other methods, the support of the estimated singular functions is smaller and in particular does not extend far beyond the area spanned by the sample data. This means that according to kernel CCA, there is large uncertainty as to which state a point in space belongs to as soon as it is outside the densely populated areas of the wedges. On the other hand, the score is lower compared to kernel EDMD or VAMPNets. This means that the metastable sets are separated less clearly, which can also be observed in the fuzziness of the transparent blue trajectory in Fig.~\ref{fig:sqrt-vamp}k and therefore the slow dynamics of the system are not represented as well as they are represented with, e.g., kernel EDMD.
		\item[(f)] \textit{VAMPNets}. As an architecture for the lobe $\bm{\chi}$ we choose a multilayer perceptron of depth $d=5$ with a rectified linear unit (ReLU) nonlinearities and $15$, $10$, $10$, $5$, and $1$ neurons, respectively. The network is trained using the Adam optimizer~\cite{kingma2014adam} with a learning rate of $10^{-3}$. We obtain a decision landscape that resembles the one of the backtransform with a perfect state separation. Also the idealized VAMP-$2$ score $s_\text{lim}$ based on the hidden transition matrix is within the standard deviation of the VAMPNet VAMP-$2$ score. The hyperparameters were chosen heuristically so that training was stable and yielded high scores.
	\end{enumerate}
	
	For the optimization of the parameters of kernel EDMD we found it crucial to first whiten the data by removing the empirical mean $\bm{\mu}$ and transforming it into the PCA basis via
	\begin{align}\label{eq:whiten}
		\mathbf{x}_t \mapsto C^{-\frac{1}{2}}(\mathbf{x}_t - \bm{\mu}),
	\end{align}
	where $C$ denotes the covariance matrix over the trajectory. The other methods were numerically more stable and applicable directly to the raw data. Whether whitening is required does not only depend on the method but in particular also on the chosen ansatz.
	
	\subsubsection{Coherent set detection}\label{sec:bickley}
	Here, we illustrate how the introduced decomposition methods can be used to detect coherent sets; i.e., sets of particles which are geometrically consistent under a forward-backward dynamic and small perturbations~\cite{froyland2010coherent,banisch2017understanding}.
	Following Ref.~\citenum{banisch2017understanding}, one can quantitatively describe coherent sets $A\subset \Omega$ under the transfer operator $\mathcal{T}$ (see, e.g.,~(\ref{eq:transfer-op-inhomogeneous})) as a set which is difficult to leave, i.e.,
	\begin{align}\label{eq:coherence}
	    \left\langle \mathcal{T}^*\mathcal{T} \frac{1_A}{\mu_s(A)}, 1_A \right\rangle_{\mu_s} \approx 1,
	\end{align}
	the probability of staying within $A$ under the forward-backward dynamic~\cite{banisch2017understanding} should be close to $1$. Here, $\mu_s(A)$ refers to the evaluation of the measure induced by the initial distribution.
	
	While dominant eigenfunctions of methods assuming time-homogeneous dynamics (cf. Fig.~\ref{fig:koopman-methods}) can be related to metastable sets, methods that may also be applied to time-inhomogeneous dynamics yield coherent sets~\cite{koltai2018optimal}. In particular, metastable sets can be understood as a special case of coherent sets. 
	
	Practically, these sets can be obtained by projecting Lagrangian data into the dominant (with respect to magnitude of singular values) left singular function space of an approximated Perron--Frobenius or Koopman operator~\cite{froyland2010coherent,klus2019kernel}, as these singular functions correspond to eigenfunctions of backward-forward dynamic $\mathcal{TT^*}$ or forward-backward dynamic  $\mathcal{T^*T}$~\cite{froyland2010coherent,klus2019kernel,wu2020variational} and therefore are an important ingredient for characterizing coherence~(\ref{eq:coherence}). Spatial proximity in the singular function space indicates membership of the same coherent set.
	
	As an application example we choose the Bickley jet, an idealized and periodically perturbed approximation of stratospheric flow which is described by a deterministic but non-autonomous system of ODEs~\cite{bickley1937lxxiii,rypina2007lagrangian}. The ODEs act on particles $\mathbf{x} = (x, y)\in \Omega = [0, 20]\times [-4, 4]$ and are given by
	\begin{align}\label{eq:bickley-odes}
		\begin{pmatrix}
			\dot x \\ \dot y
		\end{pmatrix} = \begin{pmatrix}
		-\frac{\partial\Psi}{\partial y}\\ \frac{\partial\Psi}{\partial x}
		\end{pmatrix}
	\end{align}
	with stream function
	\begin{align*}
		\Psi(x, y, t) &= c_3y - U_0L\tanh(y/L) + A_3U_0L \sech^2(y/L)\cos (k_1x)\\
		& + A_2U_0L\sech^2(y/L)\cos (k_2x - \sigma_2 t)\\
		& + A_1U_0L\sech^2(y/L)\cos (k_1x - \sigma_1t)
	\end{align*}
	and parameters chosen as in~\cite{banisch2017understanding}. The domain $\Omega$ is quasi-periodic in $x$-direction. The Bickley jet is widely used as a benchmark problem in the coherent set literature, e.g., in Refs.~\citenum{froyland2010coherent,froyland2012finite,hadjighasem2016spectral,banisch2017understanding,klus2019kernel,husic2019simultaneous}. 
	
\begin{figure}
		\centering
		\includegraphics{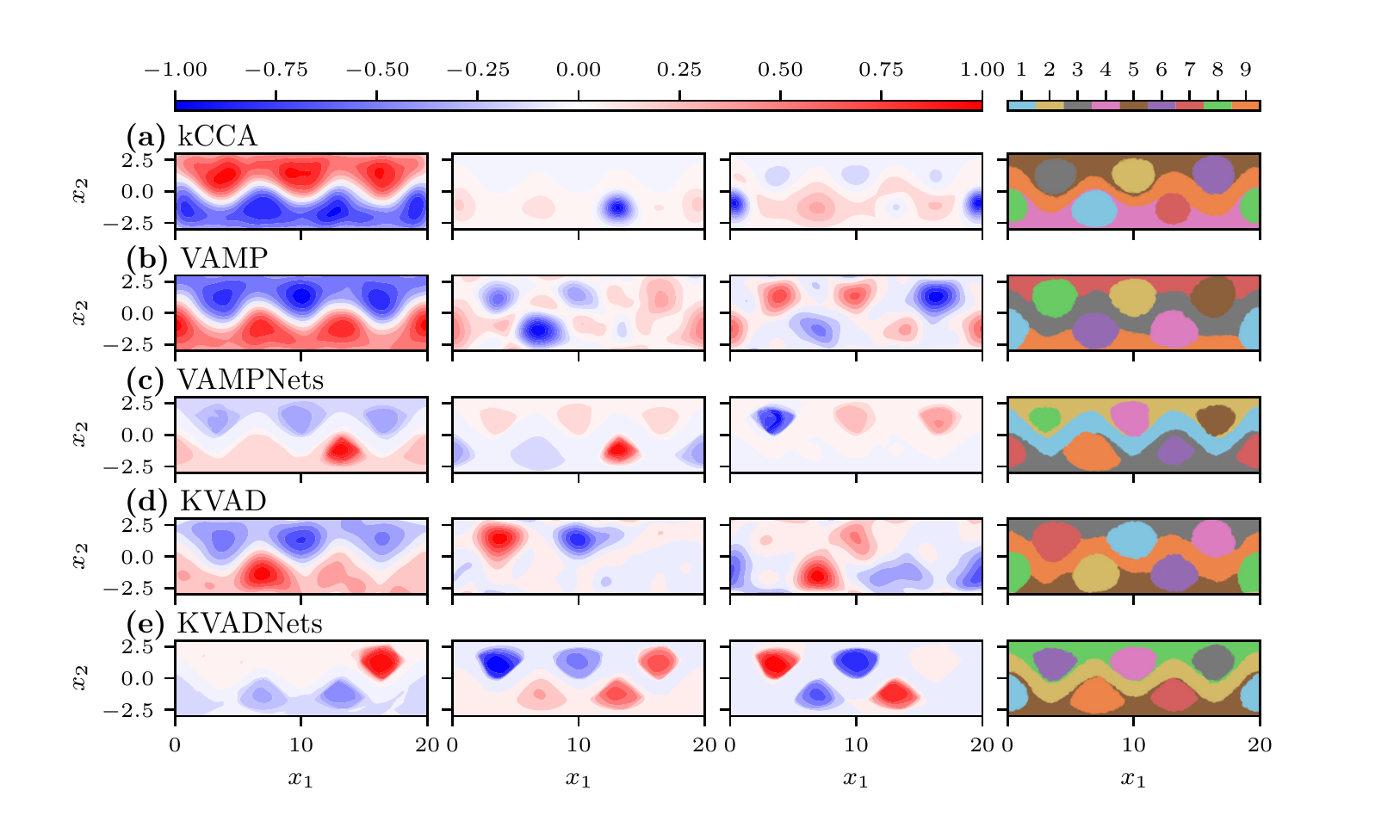}
		\caption{\textbf{Comparison of different Koopman operator methods for coherent set detection.} Columns 1--3 show the first, second, and third dominant singular function of the respective estimated Koopman operators. Column 4 shows a $k$-means clustering with $k=9$ on the initial data after it has been transformed by the first nine respective singular functions. \textbf{(a)} Kernel CCA with a Gaussian kernel with bandwidth and regularization parameter optimized to maximize the VAMP-2 score. \textbf{(b)} VAMP-estimated model, where the ansatz featurization consists of a set of randomly shifted and distorted Gaussian functions $f(x) = \exp (-x^2)$. \textbf{(c)} VAMPNets with multilayer perceptron lobes. \textbf{(d)} KVAD with the same ansatz as VAMP and a Gaussian kernel with bandwidth $\sigma = 1$. \textbf{(e)} KVADNets with a Gaussian kernel with bandwidth $\sigma = 0.5$ and a feature transformation given by multilayer perceptrons.}
		\label{fig:bickley}
	\end{figure}

	We expect to find a separation into nine coherent sets, where the domain $\Omega$ is separated into an upper $\Omega_\text{up}$ and lower $\Omega_\text{low}$ part with three circular coherent sets each, a coherent layer that is between $\Omega_\text{up}$ and $\Omega_\text{down}$ and the remainder of $\Omega_\text{up}$ and $\Omega_\text{low}$ sans the circular coherent sets, as illustrated in any of the panels of Fig.~\ref{fig:bickley} column $4$.

	Because the ODE is not autonomous~(meaning $\dot{\mathbf{x}} = \mathbf{f}(t, \mathbf{x})$ depends on time $t$), we restrict ourselves to methods that support time-inhomogeneous dynamics, in particular kernel CCA, VAMP, VAMPNets, KVAD, and KVADNets. In order to fit the respective Koopman models, we first integrate $N=3000$ particles whose positions are drawn uniformly in $\Omega$ from $t_0=0$ to $t_1=40$. From the resulting trajectories we use the initial time particle position matrix $X\in\Omega^N$ and final time particle position matrix $Y\in\Omega^N$ to find an embedding with corresponding Koopman or PF operator that describes transport from $\mathbf{x}_i$ to $\mathbf{y}_i$.
	
	The visualization of the first three estimated dominant singular functions already reveals some of the coherent structure of the underlying process~(see Fig.~\ref{fig:bickley} columns 1--3). All of the methods yield similar results and, with different degrees of sharpness, each show the three vortices in the upper part and lower part of the domain.
	
	To obtain crisp assignments to for a predefined number of coherent sets we perform $k$-means clustering using kmeans++ initialization with $k=9$ cluster centers with one cluster center belonging to exactly one coherent set. The clustering is repeated $500$ times and we select the cluster centers which yield the smallest cumulative squared distance between sample points and assigned cluster center (sometimes referred to as ``inertia''). In the last column of Fig.~\ref{fig:bickley}, the particle positions at $t=0$ are color-coded according to their cluster membership.
	
	For both VAMP and KVAD we set up the feature functions in the following way: Weight matrices $W_1\in\mathds{R}^{100\times 3}$, $W_2\in\mathds{R}^{50, 100}$ are generated by drawing i.i.d. samples from the normal distribution $\mathcal{N}(0, 1)$ and bias vectors $b_1\in\mathds{R}^{100}$, $b_1\in\mathds{R}^{50}$ are generated by drawing i.i.d. samples from the uniform distribution $\mathcal{U}(-1, 1)$. The vector-valued feature function is then given by
	\begin{align*}
	    F : \mathds{R}^2\to\mathds{R}^{50},\quad \mathbf{x}\mapsto W_2 \sigma(W_1T(\mathbf{x}) + b_1) + b_2,
	\end{align*}
	where $\sigma(x) = \exp (-x^2)$ acts component-wise and $T$ is a transformation that embeds the two-dimensional data into three dimensions by mapping it onto a cylinder
	\begin{align}\label{eq:bickley-to3d}
		T : \mathds{R}^2\to\mathds{R}^3,\quad\begin{pmatrix}
			x \\ y
		\end{pmatrix} \mapsto \begin{pmatrix}
			\cos(2\pi x / 20)\\
			\sin(2\pi x / 20)\\
			y/3
		\end{pmatrix},
	\end{align}
	accounting for the quasi-periodicity of the domain $\Omega$. KVAD is equipped with a Gaussian kernel with bandwidth $\sigma=1$.
	
	For VAMPNets, the instantaneous and time-lagged lobes are multilayer perceptrons (MLPs) and are using shared weights. The two-dimensional data is first transformed into three dimensions to account for quasi-periodicity in $x$ direction via~(\ref{eq:bickley-to3d}) and subsequently transformed through a batch normalization layer. The MLPs possess layers with $256$, $512$, $128$, $128$, and $9$ neurons, respectively, separated using ELU nonlinearities and dropout ($p=50\%$).
	
	In the case of KVADNets there is per definition just one lobe. Its architecture is the same as for VAMPNets.
	
	All models project onto the dominant nine singular functions.
	\begin{enumerate}
	    \item[(a)] \textit{Kernel CCA}. We use a Gaussian kernel where the bandwidth and regularization parameter maximize the VAMP-$2$ score ($\sigma\approx 0.58$, $\varepsilon\approx 5.6\cdot 10^{-3}$). Optimization was performed with SLSQP~\cite{kraft1988software}. The first singular function shows a clear separation between upper and lower part of the domain as the dominant process. The vortices are circular in shape and can be observed in the evaluation of the singular functions as well as the clustering.
	    \item[(b)] \textit{VAMP}. The results are qualitatively comparable to the ones of kernel CCA, however the clustering is less pronounced.
	    \item[(c)] \textit{VAMPNets}. Some of the coherent structures are clearly visible in the first three singular functions. The clustering is pronounced; however, it yields vortices of varying sizes.
	    \item[(d)] \textit{KVAD}. We use a Gaussian kernel with bandwidth $\sigma=1$. The results are qualitatively comparable to VAMP.
	    \item[(e)] \textit{KVADNets}. Here the vortices can easily be detected in the singular functions; however, the shape is less circular compared to the other methods. Furthermore, the first singular function has a less pronounced decision surface between upper and lower part of the domain; rather, it almost exclusively describes the exchange of mass between two individual vortices. The clustering, however, is sharp and is comparable to the other methods in terms of the detected sets.
	\end{enumerate}

	While differences in estimated coherent sets can be evaluated qualitatively by visual inspection, we now seek to compare the methods in a quantitative fashion and try to determine a ``best'' subdivision into coherent sets according to some criterium.
	
\begin{figure}
		\centering
		\includegraphics{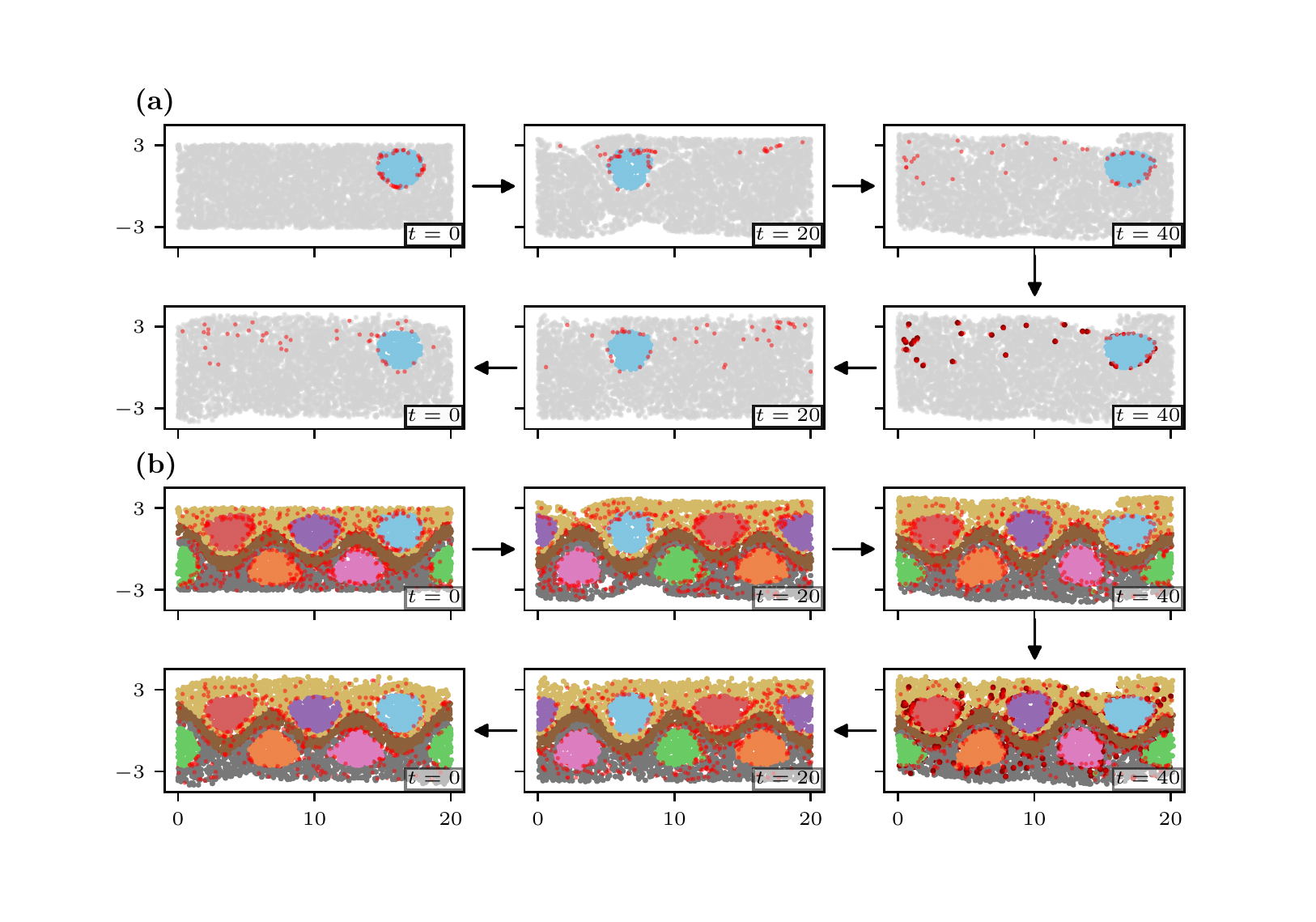}
		\caption{\textbf{Coherence in the Bickley jet.} This figure demonstrates the concept of ``leaked'' particles used for defining a coherence score. Particles are colored according to their coherent set membership or marked as leaked (red dots).
		\textbf{(a)} Particles that stay assigned to a particular coherent set (blue dots) and particles that change their set membership under the forward-backward transformation with noise (red dots) as detected by KVADNets. The top row depicts application of forward dynamics, bottom row application of noise and backward dynamics, respectively. Particles initially not belonging to the considered coherent set are shown in grey. \textbf{(b)} Generalization of (a) depicting all coherent sets; red dots now depict particles leaking from \textit{any} coherent set.}
		\label{fig:bickley_backmapped}
	\end{figure}

	To this end, we define a ``coherence score''. Let $\mathbf{x}_t = \bm{\Phi}_{t_0, t}(\mathbf{x}_0)$ be the flow describing the solution of the governing equations given an initial position $\mathbf{x}_0$ at initial time $t_0$.
	Since in this example the ground truth dynamics are known, we can take our definition of a coherent set~(\ref{eq:coherence}) as template. For a subdivision of $\Omega$ into disjoint coherent sets $\bigcup_i A_i = \Omega$, the score restricted to one set $A_i$ is defined as
	\begin{align}\label{eq:coherence-score-local}
	    s^{(i)}_{\mathrm{coh}} := \mathds{P}\left[ (\bm{\Phi}_{t_0, t_1}^{-1} \circ \mathbf{N}_\sigma \circ \bm{\Phi}_{t_0, t_1})(\mathbf{x}_{t_0}) \in A_i \mid \mathbf{x}_{t_0}\in A_i \right],
	\end{align}
	where $\mathbf{N}_\sigma(\mathbf{x}) := \mathbf{x} + \sigma\bm{\eta}$ distorts the forward-mapped $\mathbf{x}_0$ by white noise $\bm{\eta}=(\eta_1,\ldots,\eta_d)^\top,\;\eta_i\sim\mathcal{N}(0, 1)$ with standard deviation $\sigma$. In this example, we chose $\sigma=10^{-1}$. In other words, Equation~(\ref{eq:coherence-score-local}) describes the probability of a particle staying inside set $A_i\subset\Omega$ under propagation forward in time, addition of noise, and subsequent back-propagation to its initial time. This concept is illustrated in Fig.~\ref{fig:bickley_backmapped}a. Following the arrows in the figure, it depicts a subset $A_i$ at $t_0=0$ in blue and red, the remainder is colored in light gray. The particles are then propagated according to the flow $\bm{\Phi}_{t_0, t_1}$ to the final time $t_1=40$ (upper right panel). The map $N_\sigma$ is applied, yielding slightly different particle positions at $t_1$ (lower right panel). Subsequently, the distorted particles are mapped back to $t_0=0$. As one might expect, the particles leaving $A_i$ aggregate at the set's boundary. The figure uses the subdivision of the domain that is yielded by the KVADNets trained transfer operator~(see Fig.~\ref{fig:bickley}e).

	Finally, Fig.~\ref{fig:bickley_backmapped}b shows the forward-backward mapping for each of the coherent sets. For clarity, we do not distinguish leaked particles from the respective sets but only show them as generally leaked from any of the sets. It can be observed that most of the interior area of the detected vortices remains vacant of leaked particles.
	
	In order to arrive at one value for all estimated coherent sets, we consider the expectation
	\begin{align}\label{eq:coherence-score}
	    s_\mathrm{coh} := \mathds{E}_{\mu_{t_0}}\left[s_{\mathrm{coh}}^{(i)}\right] = \sum_i\frac{\mu_{t_0}(A_i)}{\mu_{t_0}(\Omega)} s_{\mathrm{coh}}^{(i)}.
	\end{align}
	For practical evaluation of the score, we estimate a MSM without a reversibility constraint (see Section~\ref{sec:msm}) on $n=2500$ discrete trajectories with nine states corresponding to the coherent sets, each trajectory corresponding to one individual particle and containing exactly two entries; namely, the coherent set before and after application of the forward-backward dynamic as given in~(\ref{eq:coherence-score-local}). Then, the $i$th diagonal entry of the transition matrix is exactly the coherence score~(\ref{eq:coherence-score-local}) corresponding to the $i$th coherent set $A_i$. The distributions $\mu_{t_0}(A_i)$ and $\mu_{t_0}(\Omega)$ are the empirical distributions; i.e., the number of particles initially inside set $A_i$ and the total number of particles $N$, respectively.
	
	A drawback of this score is that it does not indicate how faithfully the discovered coherent sets represent the system's dynamics. If, for example, the entire domain is its own coherent set, the score is maximized. Therefore, the score~(\ref{eq:coherence-score}) should be considered in conjunction with other indicators such as the VAMP or KVAD score.
	
	We compare these metrics in Table~\ref{tab:bickley-scores}.  For all used methods of this example it shows the coherence score, the VAMP-$2$ score, and the KVAD score calculated with a Gaussian kernel with bandwidth $\sigma=0.5$.

	\begin{table}	
		\begin{center}
			\resizebox{\columnwidth}{!}{
				\begin{tabular}{llllll}
					\toprule[1.5pt]
					& KVAD            & VAMP            & VAMPNets        & Kernel CCA      & KVADNets        \\
					\midrule
					Coherence score & $0.74 \pm 0.01$ & $0.77 \pm 0.01$ & $0.79 \pm 0.01$ & $0.85 \pm 0.01$ & $0.87 \pm 0.01$ \\
					VAMP-2 score  & $4.63 \pm 0.06$ & $5.18 \pm 0.08$ & $7.28 \pm 0.06$ & $5.77 \pm 0.08$ & $6.03 \pm 0.09$ \\
					KVAD score  &  $0.070 \pm 1.2\cdot 10^{-3}$ & $0.073 \pm 1.1\cdot 10^{-3}$ & $0.078 \pm 1.1\cdot 10^{-3}$ & $0.080 \pm 1.4\cdot 10^{-3}$ & $0.087 \pm 1.2\cdot 10^{-3}$ \\
					\bottomrule[1.5pt]
				\end{tabular}
			}
		\end{center}
		\caption{Table showing different kinds of scores with standard deviations for different methods used for coherent set detection using the Bickley jet example (Section~\ref{sec:bickley}). For the evaluation of the KVAD score, a Gaussian kernel with bandwidth $\sigma=0.5$ was chosen. The methods are in ascending order from left to right according to their coherence score.}
		\label{tab:bickley-scores}
	\end{table}
	
	The table also reports standard deviations, which are obtained by repeating the scoring for fifteen rounds of $n=2500$ independently and over the domain uniformly sampled initial particle positions. According to the coherence score, KVADNets deliver the best subdivision into coherent sets, with kernel CCA as a close second. This is also reflected in their respective VAMP-2 and KVAD scores. In contrast to KVAD, which yields the same sequence of methods~(if ordered ascendingly) as the coherence score, the VAMP-2 score for VAMPNets is an outlier. The VAMP-2 score for VAMPNets is significantly higher compared to any of the other methods. The Bickley jet is a deterministic system and therefore the Koopman operator associated to it is not Hilbert--Schmidt---a violation of the assumptions that are made to define the VAMP scores. Up to noise effects caused by numerical integration, this might be the cause of the high score of VAMPNets.

	It should be noted that the coherence score can still be approximated if the ground truth is not known or too expensive to compute by using a propagation model of the form~(\ref{eq:koopman-matrix}). Assuming a good representation of the slow dynamics~(which is indicated by a high VAMP or KVAD score), the error of integrating backwards in time with~(\ref{eq:koopman-matrix}) is small.

	\section{Markov state models}\label{sec:msm}
	
	Markov state models (MSMs) are stochastic models describing the time evolution of a random process $\{\mathbf{x}_t\}_{t\geq 0}$, $\mathbf{x}_t\in\Omega$~(see Refs.~\citenum{schutte_direct_1999,swope_describing_2004,singhal_using_2004,noe_hierarchical_2007,noe_probability_2008,noe_constructing_2009,prinz_markov_2011,chodera_markov_2014,husic_markov_2018}) and describe Markov chains with memory depth of 1.
	In other words, given a sequence $(..., \mathbf{x}_{t-2\tau}, \mathbf{x}_{t-\tau}, \mathbf{x}_t)$ with a set of possible states $\Omega$, the conditional probability of encountering a particular state $\mathbf{x}_{t+\tau} \in \Omega$ is only conditional on $\mathbf{x}_t \in S$; i.e. $\mathds{P}(\mathbf{x}_{t+\tau}|\mathbf{x}_t, \mathbf{x}_{t-\tau}, \mathbf{x}_{t-2\tau}, ...) = \mathds{P}(\mathbf{x}_{t+\tau}|\mathbf{x}_t)$. In contrast to the methods presented in Section~\ref{sec:dimredux}, we assume that we have a finite number of discrete states. Therefore we consider
	\begin{align}\label{eq:msm-state-set}
	    S := \{1,\ldots, n\} \cong \Omega
	\end{align}
	as state space for the remainder of this section. Often we are presented with data that does not live in a countable or even finite state space. In these cases, the state space is tessellated using a finite amount of indicator functions. Typically, the tesselation is a Voronoi decomposition. 
	
	\begin{figure}
		\centering  
		\includegraphics[scale=1.2]{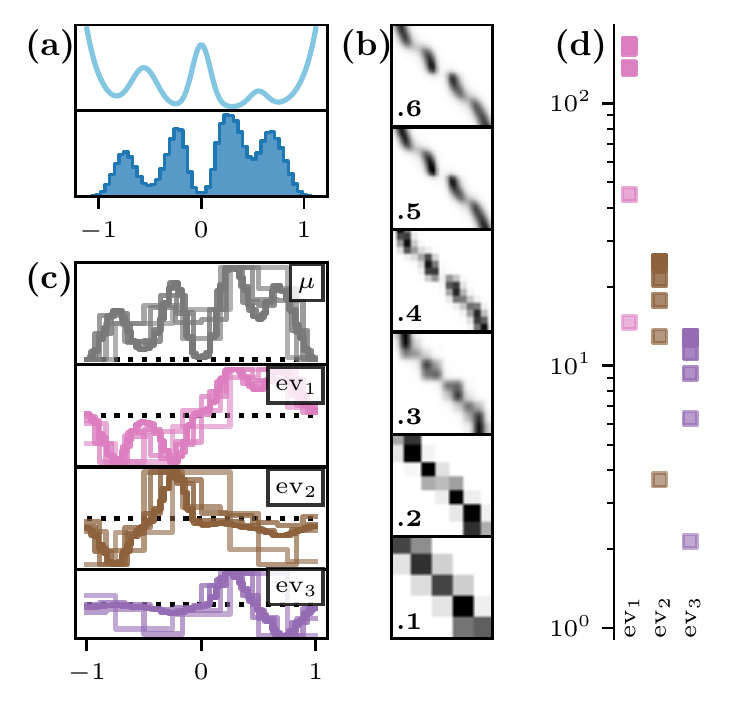}
		\caption{MSM spectral decomposition for a random walk in a asymmetric 1D 4-well potential; the corresponding potential function~(upper part) with histogram of simulated data~(lower part) is depicted in \textbf{(a)}. \textbf{(b)} shows transition matrix estimates for various discretizations, from very coarse \textbf{(b)}.1 to very fine \textbf{(b)}.6. The discrete states are sorted in ascending order with respect to their corresponding $x$-coordinate. \textbf{(c)} depicts the four dominant left eigenvectors; discretizations are color coded from faint (coarse) to strong colors (fine discretization). The eigenvalues corresponding to the nonstationary eigenvectors are depicted in \textbf{(d)} in the same color.}
		\label{fig:msm_example}
	\end{figure}
	
	As shown at the example of the Prinz potential (Fig.~\ref{fig:msm_example}a), the fineness of the chosen discretization affects MSM approximation quality~\cite{prinz_markov_2011}. The estimated transition matrix can approximate the dynamics with higher spatial resolution in a finer discretized space (Fig.~\ref{fig:msm_example}b).
	Furthermore, the discretization has implications on the estimated eigenfunctions and, in particular, the estimated stationary distribution (Fig.~\ref{fig:msm_example}c).
	Evaluating the eigenvalues for a given discretization yields a comprehensive picture of the model's quality (Fig.~\ref{fig:msm_example}d) as the true eigenvalues present an upper bound to the estimated ones (variational principle~\cite{wu2020variational}):~the sum of eigenvalues reflects the VAMP-$1$ score.

	MSMs fit into the framework of transfer operators as introduced in Section~\ref{sec:dimredux}~(see Fig.~\ref{fig:koopman-methods}). In particular, an indicator function ansatz used with VAC and/or EDMD yields an MSM. When indicator functions are used with VAMP, one obtains generalized MSMs (GMSMs), which are capable of representing time-inhomogeneous dynamics. We suggest Refs.~\citenum{husic_markov_2018, chodera_markov_2014, prinz_markov_2011} for thorough reviews.
	
	The conditional probabilities in the MSM framework are described by a transition matrix $P \in \mathds{R}^{n\times n}$, where $n=|S|$ is the number of states. The transition matrix is given by
	\begin{align}\label{eq:msm-transition-matrix}
		P_{ij} = \mathds{P}(x_{t+\tau} = j | x_t = i) \qquad \forall t\geq 0,
	\end{align}
	i.e., the time-stationary probability of transitioning from state $i\in S$ to state $j\in S$ within time $\tau$. This also means that $P$ is a row-stochastic matrix. 
	Note that the MSM transition matrix is a special case of a transfer operator approximation (see Section~\ref{sec:dimredux}), where the ansatz consists of indicator functions.
	
	Dynamical quantities of interest can be computed from an MSM's transition matrix, e.g., mean first passage times and fluxes among (sets of) states~\cite{metzner_transition_2009}, implied timescales~\cite{prinz_markov_2011}, or metastable decompositions of Markov states~\cite{roblitz_fuzzy_2013}.
	
	\paragraph{MSM estimation with deeptime}
	The goal of the \mintinline{python}{deeptime.markov} module is to provide tools to estimate and analyze MSMs from discrete-state time-series data. If the data's domain is not discrete, classical discretization algorithms (such as the ones implemented in \mintinline{python}{deeptime.clustering}) can be employed to assign each frame to a state.

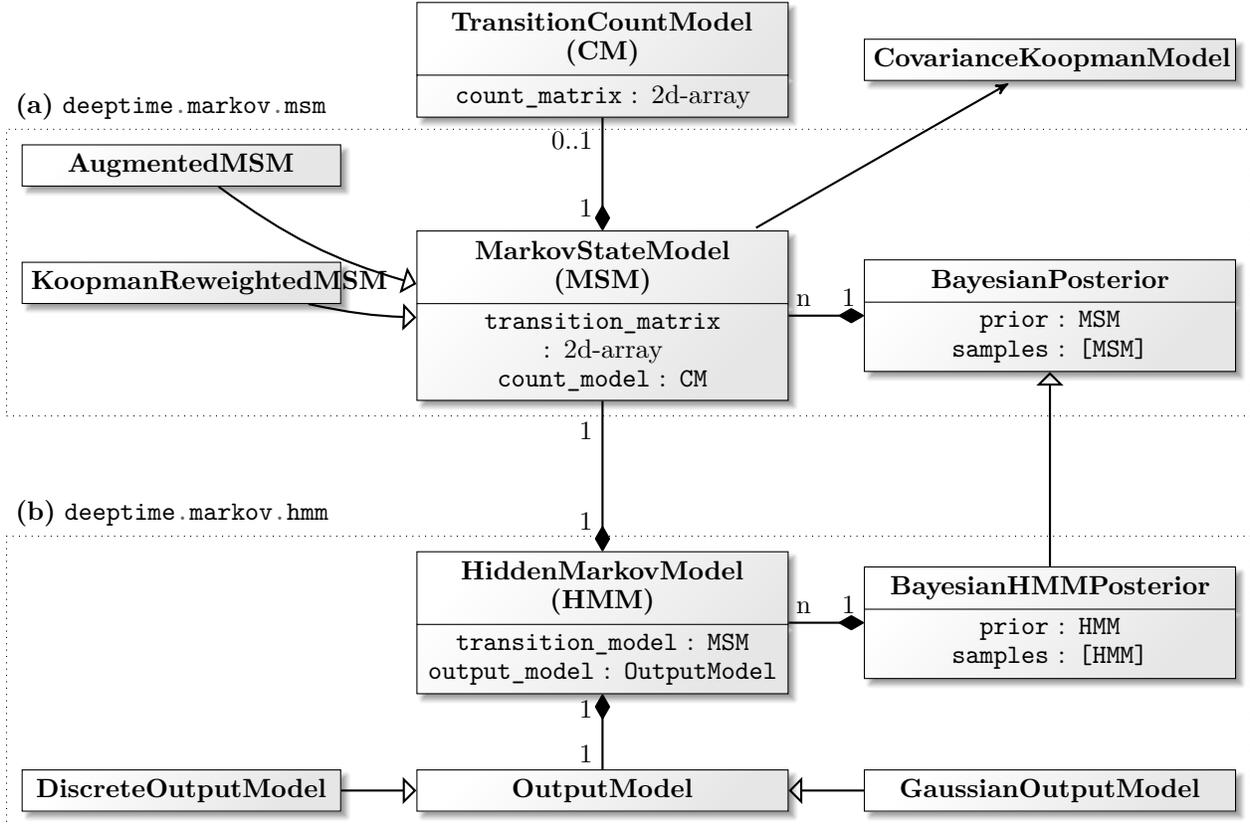
\begin{figure}
		\begin{center}
			\begin{tikzpicture}
				\node (transitions) [class2, rectangle split, rectangle split parts=2]
				{
					\textbf{TransitionCountModel (CM)}
					\nodepart{second}\mintinline{python}{count_matrix} : 2d-array
				};
				\node(aux) [text width=4cm, left=of transitions] {};
				\node (msm) [class2, below=1.5cm of transitions, rectangle split, rectangle split parts=2]
				{
					\textbf{MarkovStateModel (MSM)}
					\nodepart{second}
					\mintinline{python}{transition_matrix} : 2d-array\\
					\mintinline{python}{count_model} : \mintinline{python}{CM}
				};
				\node (hmm) [class2, below=2cm of msm, rectangle split, rectangle split parts=2]
				{
					\textbf{HiddenMarkovModel (HMM)}
					\nodepart{second}
					\mintinline{python}{transition_model} : \mintinline{python}{MSM}\\
					\mintinline{python}{output_model} : \mintinline{python}{OutputModel}
				};
				\node (outputmodel) [class2, below=of hmm]
				{
					\textbf{OutputModel}
				};
				\node (discreteoutputmodel) [class2, left=of outputmodel, text width=4cm]
				{
					\textbf{DiscreteOutputModel}
				};
				\node (gaussianoutputmodel) [class2, right=of outputmodel]
				{
					\textbf{GaussianOutputModel}
				};
				\node (covkoopman) [class2, right=of transitions]
				{
					\textbf{CovarianceKoopmanModel}
				};
				\node (bmsm) [class2, right=of msm, rectangle split, rectangle split parts=2]
				{
					\textbf{BayesianPosterior}
					\nodepart{second}
					\mintinline{python}{prior} : \mintinline{python}{MSM}\\
					\mintinline{python}{samples} : \mintinline{python}{[MSM]}
				};
				\node (bhmm) [class2, right=of hmm, rectangle split, rectangle split parts=2]
				{
					\textbf{BayesianHMMPosterior}
					\nodepart{second}
					\mintinline{python}{prior} : \mintinline{python}{HMM}\\
					\mintinline{python}{samples} : \mintinline{python}{[HMM]}
				};
				\node (amm) [class2, below=of aux, text width=4cm]
				{
					\textbf{AugmentedMSM}
				};
				\node (oom) [class2, below=of amm, text width=4cm]
				{
					\textbf{KoopmanReweightedMSM}
				};
				\path[comp,->] (transitions) edge node[left,pos=0.2] {0..1} node[left,pos=.8] {1} (msm);
				\path[comp,->] (msm) edge node[left,pos=0.2] {1} node[left,pos=.8] {1} (hmm);
				\path[comp,->] (msm) edge node[above,pos=0.2] {n} node[above,pos=.8] {1} (bmsm);
				\path[comp,->] (hmm) edge node[above,pos=0.2] {n} node[above,pos=.8] {1} (bhmm);
				\path[comp,->] (outputmodel) edge node[left,pos=0.2] {1} node[left,pos=.8] {1} (hmm);
				\path[inharrow,->] (bhmm) edge (bmsm);
				\path[inharrow,->] (amm) edge[bend right=10] (msm);
				\path[inharrow,->] (oom) edge[bend right=5] (msm);
				\path[inharrow,->] (discreteoutputmodel) edge (outputmodel);
				\path[inharrow,->] (gaussianoutputmodel) edge (outputmodel);
				\path[pil,->] (msm) edge (covkoopman);
				\node[draw,dotted,fit=(amm) (oom) (bmsm) (msm),inner sep=2mm] (box1) {};
				\node[above right] at (box1.north west){\textbf{(a)}~\mintinline{python}{deeptime.markov.msm}};
				\node[draw,dotted,fit=(hmm) (bhmm) (discreteoutputmodel) (gaussianoutputmodel),inner sep=2mm] (box2) {};
				\node[above right] at (box2.north west){\textbf{(b)}~\mintinline{python}{deeptime.markov.hmm}};
			\end{tikzpicture}
		\end{center}
		\caption{\textbf{Class diagram showing relationships of main models contained in the \mintinline{python}{markov} package.} We show two of the sub-packages of the \mintinline{python}{markov} module. Classes may be related via inheritance ({\protect\tikz[baseline=-0.5ex]\protect\draw[inharrow] (0mm,0mm) -- +(4mm,0mm);}), composition ({\protect\tikz[baseline=-0.5ex]\protect\draw[comp] (0mm,0mm) -- +(4mm,0mm);}), or produce objects of another kind ({\protect\tikz[baseline=-0.5ex]\protect\draw[pil] (0mm,0mm) -- +(4mm,0mm);}). In case of composition, we further denote the cardinality of the relationship next to the arrow.
			\textbf{(a)} \mintinline{python}{msm}. The \mintinline{python}{markov.msm} package has the Markov state model (MSM) at its core. It is completely determined by a transition matrix, but may also contain information about the statistics of data, in which case it possesses a \mintinline{python}{TransitionCountModel}. Furthermore there are \mintinline{python}{AugmentedMSM} and \mintinline{python}{KoopmanReweightedMSM} subclasses. An MSM is also a Koopman model using indicator basis functions; therefore, it can generate corresponding \mintinline{python}{CovarianceKoopmanModel} objects (see Section~\ref{sec:dimredux} and Fig.~\ref{fig:koopman-inheritance}). Bayesian sampling around an MSM leads to \mintinline{python}{BayesianPosterior}s, consisting of the prior and drawn samples. \textbf{(b)} \mintinline{python}{hmm}. This package contains estimators and models corresponding to hidden Markov model estimation. A \mintinline{python}{HiddenMarkovModel} consists out of a \mintinline{python}{transition_model} which describes evolution of a hidden state and an \mintinline{python}{OutputModel} which assigns a distribution over observable states to each hidden state. The discrete output model assigns each hidden state a discrete probability distribution over states, the Gaussian output model samples from one-dimensional Gaussians with means and variances conditioned on the hidden state. Same as MSMs, Bayesian sampling is available for HMMs, leading to a \mintinline{python}{BayesianHMMPosterior}.}
		\label{fig:markov-models}
	\end{figure}

	In what follows, we introduce the core object, the \mintinline{python}{MarkovStateModel}, as well as a variety of estimators. An overview of the main models contained in the \mintinline{python}{markov} module is depicted in Fig.~\ref{fig:markov-models}.
	
	Deeptime implements maximum-likelihood estimators for Markov state models as well as Bayesian sampling routines~\cite{trendelkamp-schroer_estimation_2015}, leading to \mintinline{python}{MarkovStateModel} and \mintinline{python}{BayesianPosterior} model instances, respectively. An integral component of MSM estimation and sampling based on time-series data is collecting statistics over the encountered state transitions (transition counting), which leads to a \mintinline{python}{TransitionCountModel}.  
	
	Bayesian sampling of MSMs leads to a \mintinline{python}{BayesianPosterior} that consists out of one \mintinline{python}{MarkovStateModel} instance representing the prior as well as the sampled \mintinline{python}{MarkovStateModel}s instances (see Fig.~\ref{fig:markov-models}a). Each \mintinline{python}{MarkovStateModel} possesses a transition matrix~(\ref{eq:msm-transition-matrix}) and, if available, statistical information about the data in the form of a transition count matrix. Furthermore, deeptime provides augmented Markov models (AMMs)~\cite{olsson_combining_2017} which can be used when experimental data is available, as well as observable operator model MSMs (OOMs)~\cite{nuske_markov_2017}. OOMs are unbiased estimators for the MSM transition matrix that correct for the effect of being presented with out of equilibrium data even when short lag-times are used. Both AMMs and OOMs inherit from the \mintinline{python}{MarkovStateModel} class definition.

	While MSMs are a special case of the transfer operator model (cf.~Section~\ref{sec:dimredux-shallow} and Figures~\ref{fig:koopman-methods} and~\ref{fig:koopman-inheritance}), they can be converted to \mintinline{python}{CovarianceKoopmanModel}s of two different types. 
	In one case, one can define the Koopman operator solely based on the transition matrix and corresponding stationary distribution; i.e., without respect to any statistical information. In the other case, when statistical information is present in the form of a \mintinline{python}{TransitionCountModel}, the statistics over transition counts may be used to estimate an empirical distribution according to which the Koopman operator is defined. The choice of Koopman model is up to the user; therefore, in deeptime MSMs do not inherit from \mintinline{python}{CovarianceKoopmanModel} but rather offer properties yielding respective instances of \mintinline{python}{CovarianceKoopmanModel}.
	
	When estimating MSMs from data, deeptime assumes that the data is in the form of $k\geq 1$ trajectories $T_1, T_2, \ldots, T_k$ which comprise sequences of discrete states, i.e.,
	\begin{align}
		T_i = (s_1, s_2, \ldots, s_{n_i}), \; \forall j = 1,\ldots, n_i : s_j \in S,
	\end{align}
	where $n_i$ is the length of the $i$-th trajectory and $S = \{ 0, 1, \ldots, N_S - 1\}$ is the set of discrete states. In terms of further analysis it can be desirable to restrict the discrete state space onto a subset of $S'\subset S$, e.g., when certain state transitions are not populated and/or to select an ergodic subset. This task is best performed using a \mintinline{python}{TransitionCountModel} instance prior to estimating an MSM, as it possesses methods to produce new instances of the transition count model but restricted onto $S'$.
	
	\paragraph{Hidden Markov models}
	In many applications, the observed processes are only approximately Markovian in discrete state space; i.e., MSMs are only approximately valid~\cite{prinz_markov_2011}.
	The Markovianity assumption for the observed dynamics is discarded for hidden Markov models (HMMs) which assume that the modeled stochastic process is hidden~(not directly observable).
	Therefore, the central object of the HMM is the transition matrix $\tilde P$ among hidden states $s_i \in S$.
	The transition matrix $\tilde P$ can be estimated from the time series of observable states $O$ with the Baum--Welch algorithm~\cite{baum1966statistical, baum1967inequality, baum1970maximization, baum1972inequality}. 
	Briefly:~alongside the transition matrix $\tilde P$, for each hidden state $s_t \in S$ the algorithm estimates an emission probability for a given observable state $\mathbf{o}_t \in O$. 
	HMMs therefore provide a (time-dependent) mapping between observable and hidden states along with the transition matrix $\tilde P$~\cite{rabiner_tutorial_1989}.
	This further allows us to estimate a maximum likelihood pathway of the trajectories in the hidden state space (Viterbi algorithm~\cite{viterbi_error_1967}).
	
	Because the Baum--Welch algorithm converges to a local likelihood maximum~\cite{rabiner_tutorial_1989}, it is crucial to provide a reasonable initial guess of the emission probabilities and initial state distribution. 
	Deeptime offers multiple possibilities to initialize the HMM estimation procedure (contained in the \mintinline{python}{deeptime.markov.hmm.init} package), with a fallback option to a classical MSM or MSM-derived~(e.g., PCCA~\cite{roblitz_fuzzy_2013}) estimate of the metastable dynamics.
	
	The initial guess is an object of type \mintinline{python}{HiddenMarkovModel}~(see Fig.~\ref{fig:markov-models}b). HMMs are composed of an MSM which describes the hidden state transitions and an output model. The output model is responsible for mapping a hidden state $s_t$ to an observable state $\mathbf{o}_t = \mathbf{o}(s_t)\in O$. Deeptime offers \mintinline{python}{DiscreteOutputModel}s which map each hidden state to a sample of a discrete probability distribution over observable states as well as \mintinline{python}{GaussianOutputModel}s which map a hidden state to a sample of a one-dimensional Gaussian distribution with mean and variance depending on the hidden state.
	
	As with MSMs, HMMs in deeptime also support Bayesian sampling following a Gibbs sampling scheme detailed in Ref.~\citenum{chodera2011bayesian}. This produces a \mintinline{python}{BayesianHMMPosterior} which inherits from the \mintinline{python}{BayesianPosterior},~(cf. Fig.~\ref{fig:markov-models}), which allows samples of quantities of interest which can be derived from an HMM instance to be collected.
	
	\section{Sparse identification of nonlinear dynamics}\label{sec:sindy}
	
	The sparse identification of nonlinear dynamics (SINDy) algorithm~\cite{brunton2016sindy} is a data-driven method for discovering nonlinear dynamical systems models from measurement data using sparse regression. The method also fits into the Koopman operator framework presented in Section~\ref{sec:dimredux} since it is related to an approximation of the Koopman generator, defined by $\mathcal{L}f := \lim_{\tau\to 0} (\mathcal{K}_\tau f - f) / \tau$, see Ref.~\citenum{klus2020data} for details. The goal of SINDy is to approximate a nonlinear dynamical system
	\begin{align}
		\frac{d}{dt}\mathbf{x} = \mathbf{f}(\mathbf{x})
	\end{align}
	as a sparse linear combination of candidate functions $\theta_k(\mathbf{x})$:
	\begin{align}
		\frac{d}{dt}x_j \approx \sum_{j=1}^\ell \xi_{jk}\theta_k(\mathbf{x}) \quad\Longrightarrow\quad \frac{d}{dt}\mathbf{x} \approx \Xi\mathbf{\Theta}(\mathbf{x}).
	\end{align}
	The matrix $\Xi$ is assumed to be sparse, with the nonzero elements determining which terms in the library $\mathbf{\Theta}$ are active in the dynamics.  
	In practice, the library $\mathbf{\Theta}$ is defined either to contain a generic set of terms, such as monomials, or terms guided by partial knowledge of the physical system. 
	For example, metabolic regulatory networks often include rational function nonlinearities~\cite{Mangan2016ieee}.  
	However, monomomials often suffice, either because the governing physics is polynomial (e.g., the Navier--Stokes equations for fluid dynamics), or because polynomials provide a reasonable Taylor expansion of the dynamics into a normal form~\cite{Champion2019pnas}.  
	
	The sparse matrix $\Xi$ is typically identified via sparse regression based on a time-series of data $\mathbf{x}_1, \mathbf{x}_2, \dots, \mathbf{x}_m$ collected at times $t_1, t_2, \dots, t_m$.  
	This data is organized into a matrix $X\in\mathds{R}^{n\times m}$, and a matrix of derivatives $\dot{X}$ is formed either by measuring the derivatives directly or numerically approximating them from the data in $X$. 
	The library $\mathbf{\Theta}$ may now be evaluated on the data matrix $X$, resulting in the following matrix system of equations
	\begin{align}\label{eq:sindy}
		\dot{X} \approx \Xi\mathbf{\Theta}(X).
	\end{align}
	The matrix $\Xi$ is then solved for in the following optimization
	\begin{align}
		\argmin_{\Xi} \| \dot{X} - \Xi\mathbf{\Theta}(X)\|_F + \lambda\| \Xi\|_0.
	\end{align}
	The first term measures the model error, while the $\|\cdot\|_0$ term counts the number of nonzero elements in $\Xi$, promoting sparsity. 
	This zero norm is non-convex, and several relaxations are available that yield sparse solutions~\cite{brunton2016sindy,champion2020unified}. 
	
	There are several extensions to SINDy, e.g., incorporating the effect of actuation and control~\cite{Kaiser2018prsa} and to enforce partially known physics, such as symmetries and conservation laws~\cite{Loiseau2017jfm}.  
	It is also possible to combine SINDy with deep autoencoders to identify a coordinate system in which the dynamics are approximately sparse~\cite{Champion2019pnas}.  
	Other extensions include the discovery partial differential equations~\cite{Rudy2017sciadv,Schaeffer2017prsa}, the modeling of stochastic dynamics~\cite{boninsegna2018sparse,klus2020generator,callaham2021nonlinear}, and weak formulations of the problem~\cite{Schaeffer2017pre,gurevich2019robust,Reinbold2020pre}, among others~\cite{tran2017exact,schaeffer2017learning,Schaeffer2017prsa,zhang2019convergence}. 
	SINDy has also been extended to accommodate tensor libraries, which dramatically increases its ability to handle systems with high state dimension~\cite{Gelss2019mindy}. 
	This sparse modeling procedure has been applied to discover new physical models, for example in fluid dynamics~\cite{Loiseau2017jfm,Deng2020JFM}, including for turbulence closure modeling~\cite{beetham2021sparse}.
	
	It is important to note that SINDy also applies equally well to discrete time systems
	\begin{align}
		\mathbf{x}_{k+1} = \mathbf{F}(\mathbf{x}_k)
	\end{align}
	in which case derivatives need not be estimated.  
	If SINDy is formulated in discrete time with no sparsity promoting term (i.e., $\lambda=0$) and with a library $\mathbf{\Theta}(\mathbf{x})=\mathbf{x}$, then the DMD approximation is recovered.  
	
\begin{figure}
		\centering
		\includegraphics[width=.5\textwidth]{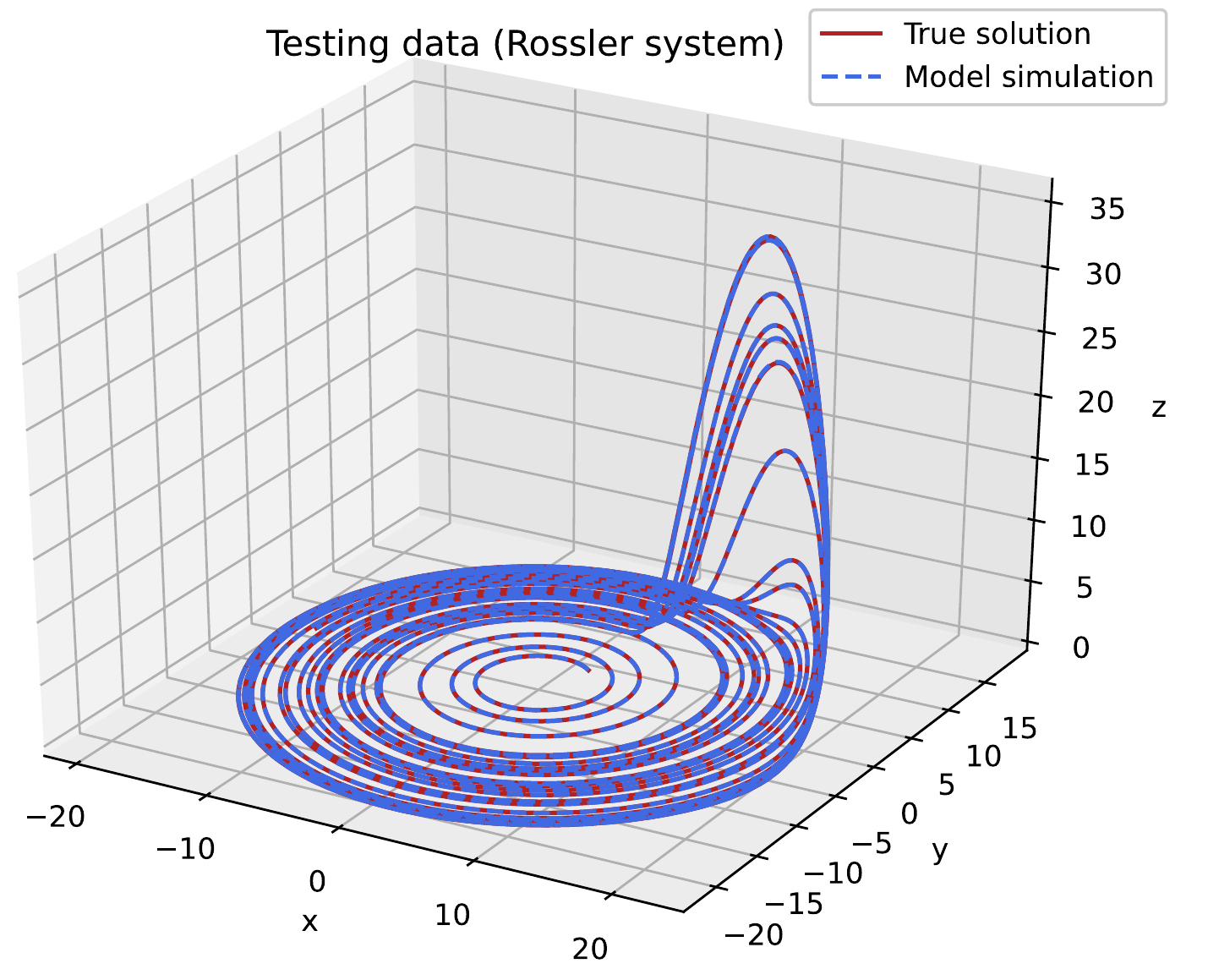}
		\caption{Reconstruction of the Rössler attractor using SINDy.}
		\label{fig:SINDyRossler}
	\end{figure}

	To demonstrate SINDy, we consider the Rössler attractor~\cite{rossler1976equation}, a system of ODEs exhibiting chaotic behavior. Fig.~\ref{fig:SINDyRossler} shows the reconstruction of the dynamic attractor for the Rössler system of equations:
	\begin{align*}
			\dot{x}_1 & = -x_2 - x_3\\
			\dot{x}_2 & = x_1 + ax_2\\
			\dot{x}_3 & = b + x_3(x_1-c)
	\end{align*}
	with constants $a=0.1$, $b=0.1$, and $c=14$.
	
	Deeptime has two SINDy objects. The \mintinline{python}{SINDy} estimator is used to solve the optimization problem~(\ref{eq:sindy}) given $\mathbf{\Theta}$, $X$, and optionally $\dot{X}$. By default the sequentially-thresholded least-squares algorithm~\cite{brunton2016sindy} is used to solve the optimization problem. If $\dot{X}$
	is not user-provided, it is estimated from $X$ with a first order finite difference method.
	
	The estimator produces a \mintinline{python}{SINDyModel}, representing the learned dynamical system. The model can be used to predict derivatives given state variables, to simulate forward in time from novel initial conditions, and to score itself against ground truth data.
	
	The implementation is API-compatible to the Python package PySINDy~\cite{desilva2020}, which in particular enables users to make use of a wider range of optimizers defined in PySINDy.

	\section{Datasets}\label{sec:data}
	
	Deeptime offers a range of datasets to which its methods can be applied. The datasets and methods were purposefully designed to be non-domain-specific and to deliver data generators rather than fixed datasets.
	As a result, the repository as well as package size are remain small and generation parameters can be varied to study their effects on the algorithms. The data simulators are structured so that performance-critical parts are implemented in C++ and the generation procedure is not very time consuming.
	
	In particular, a range of example SDEs of the form
	\begin{align*}
	    \mathrm{d}\mathbf{x}_t = \mathbf{F}(t, \mathbf{x}_t)\mathrm{d}t + \sigma\mathrm{d}W_t,
	\end{align*}
	where $\mathbf{F}:\mathds{R}\times\mathds{R}^d\to\mathds{R}^d$, $W_t$ a $d$-dimensional Wiener process, and $\sigma\in\mathds{R}^{d\times d}$, are implemented. All these SDEs are integrated using an Euler--Maruyama integrator. While the definition of these examples happens in C++, it is set up in such a way that also C++-inexperienced users can natively define their own.
	
\begin{figure}
		\centering
		\includegraphics[width=\columnwidth]{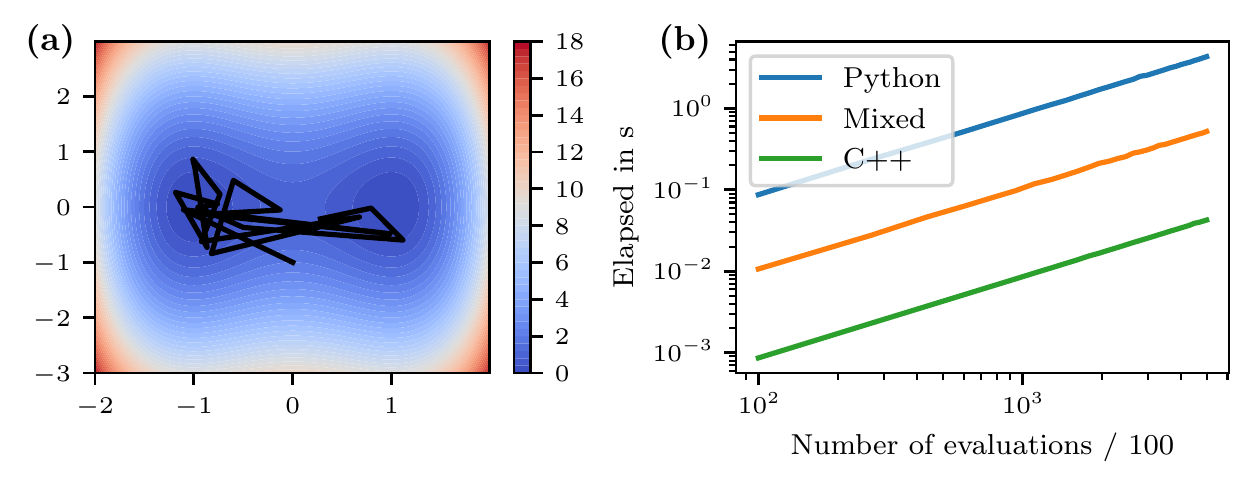}
		\caption{\textbf{Two-dimensional double-well example system.} We show a performance comparison between three different implementations of a two-dimensional double well system. \textbf{(a)} Potential landscape $V(\mathbf{x}) = V(\mathbf{x}_1, \mathbf{x}_2) = (\mathbf{x}_1^2 - 1)^2 + \mathbf{x}_2^2$ with example trajectory under diffusion matrix $\sigma=\mathrm{diag}(0.7, 0.7)$ and state snapshots taken every $10^4$ steps. \textbf{(b)}~Time elapsed over number of evaluations, while one evaluation corresponds to evolving the state by $100$ steps under a integration step size of $h=10^{-3}$. ``Python'' refers to a native Python implementation, ``C++'' refers to a native C++ implementation, and ``mixed'' refers to a C++ implementation, where the diffusion matrix as well as the gradient of the potential are defined in Python.}
		\label{fig:potential_performance}
	\end{figure}

	For example, the definition of a double well system
	\begin{align*}
	    \mathrm{d}\mathbf{x}_t = -\nabla V(\mathbf{x}_t) \mathrm{d}t + \sigma\mathrm{d}W_t,\quad V(\mathbf{x}) = (\mathbf{x}_1^2 - 1)^2 + \mathbf{x}_2^2,
	\end{align*}
	with $\mathbf{x}_t\in\mathds{R}^2$ and $\sigma = \mathrm{diag}(0.7, 0.7)$ can be achieved by a struct definition detailing the evaluation of the right-hand side. Many of the parameters of the system can be made available at compile-time, enabling further optimizations by the compiler. An example trajectory as well as a contour plot of the potential landscape can be found in Fig.~\ref{fig:potential_performance}a. By making information such as the data type (e.g., \mintinline{cpp}{float} or \mintinline{cpp}{double}), the dimension of the state space, the integrator, and $\sigma$ available at compile time, the compiler can perform further optimizations and potentially vectorizations that it otherwise could not, reducing the time it needs for evaluation.
	
	Users also have the option to define the right-hand side $\mathbf{F}(\mathbf{x}_t)$ as well as the diffusion matrix $\sigma$ in Python at some performance penalty~(see Fig.~\ref{fig:potential_performance}b. Three different implementations are compared:~one native C++ implementation, one implementation where just $\sigma$ and the right-hand side are defined in Python, and one native Python implementation. One can see that roughly one order of magnitude in terms of evaluation performance is gained from native Python to a mixed Python/C++ implementation and from the mixed implementation to a native C++ implementation.
	
	A drawback of making this information known at compile time is that for the mixed Python/C++ implementations, the dimension needs to be predefined; i.e., it must be explicitly exported when generating the Python bindings. On the other hand it improves performance and one can first prototype a system using the Python-defined diffusion matrix and right-hand side, and then eventually move the implementation to native C++ with relative ease.
	
	\section{Discussion and outlook}
	We have outlined the key components of deeptime's API and discussed the corresponding theory and methods as well as their relationships, in particular transfer operator based methods which can be used for dimension reduction, coherent set detection, analysis of kinetic quantities, and discovery of governing dynamics. These applications were each demonstrated with respective examples.

	For future development we are actively looking for contributors and want to extend the currently available library of methods and datasets. For example there is a version of VAMPNets which allows the inclusion of experimental data. The SINDy module can be extended to include neural network based estimation of dynamics. Also the HMM module can be extended to support a richer set of output models. Furthermore the inclusion of more example datasets is desirable as this enables users to test and analyze existing or new methods and draw comparisons.
	
	Finally we are planning to integrate time-series specific chunking and streaming capabilities so that methods which support online learning can more easily be used with data streams.
	
	\section*{Acknowledgements}
	We acknowledge financial support from Deutsche Forschungsgemeinschaft DFG (SFB/TRR 186, Project A12 and SFB 1114, Projects A04, B06, and C03), the European Commission (ERC CoG 772230 "ScaleCell"), and the Berlin Mathematics center MATH+ (AA1-6 and AA1-10).
	Part of this research was performed while M.H., A.M., B.E.H., S.K., H.W., N.K., S.L.B., and F.N.~were visiting the Institute for Pure and Applied Mathematics (IPAM), which is supported by the National Science Foundation (Grant No. DMS-1440415).
	F.N. acknowledges the German Ministry for Education and Research (Project BIFOLD - Berlin Institute for the Foundations of Learning and Data).
	S.L.B. and N.K. acknowledge support from the National Science Foundation AI Institute in Dynamic Systems (Grant No. 2112085). 
	B.E.H.~acknowledges the Lews-Sigler Institute at Princeton University, the Princeton Center for the Physics of Biological Function, and the Princeton Center for Theoretical Science.
	H.W.~acknowledges the NSF of China (Grant No.~12171367), the Shanghai Municipal Science and Technology Commission (No.~20JC1413500), the Shanghai Municipal Science and Technology Major Project (Grant No.~2021SHZDZX0100) and the fundamental research funds for the central universities of China (Grant No.~22120210133).
	\clearpage
	\printbibliography[title={Bibliography}]
\end{document}